\documentclass{amsart}

\usepackage{amssymb}
\usepackage{epsfig}
\usepackage{psfrag}
\newtheorem{theorem}{Theorem}[section]
\newtheorem{lemma}[theorem]{Lemma}

\newtheorem{corollary}[theorem]{Corollary}

\theoremstyle{definition}

{\bf }

\theoremstyle{remark}
\newtheorem{remark}[theorem]{Remark}
\numberwithin{equation}{section}
\newcommand{\abs}[1]{\lvert#1\rvert}

\newcommand{\ul}[1]{\underline{#1}}
\newcommand{\ol}[1]{\overline{#1}}
\newcommand{\C}{\mathbb C}
\newcommand{\E}{\mathbb E}
\newcommand{\Out}{\hbox{Out\,}}
\newcommand{\Int}{\hbox{Int\,}}
\newcommand{\Aut}{\hbox{Aut\,}}
\newcommand{\Ob}{\hbox{Ob}}
\newcommand{\Ad}{\hbox{Ad}}
\newcommand{\comment}[1]{}

\begin{document}
\title{The planar algebra of diagonal subfactors}

\author{Dietmar Bisch}
\address{Vanderbilt University, Department of Mathematics, SC 1326,
Nashville, TN 37240, USA}
%
\email{dietmar.bisch@vanderbilt.edu, paramita.das@vanderbilt.edu}
\email{shamindra.k.ghosh@vanderbilt.edu}
\thanks{The authors were supported by NSF under Grant No. DMS-0301173 
and DMS-0653717.}
%
\author{Paramita Das}
%
\author{Shamindra Kumar Ghosh}

%


\keywords{G-kernel, group cocycle, planar algebra, planar operad, subfactor}

\begin{abstract}
There is a natural construction which associates to a finitely
generated, countable, discrete group $G$ and a 3-cocycle $\omega$
of $G$ an
inclusion of II$_1$ factors, the so-called
diagonal subfactors (with cocycle). In the case when the cocycle
is trivial these subfactors are well studied and their
standard invariant (or planar algebra) is known. We give
a description of the planar algebra
of these subfactors when a cocycle is present. The action of
Jones' planar operad involves the 3-cocycle $\omega$ explicitly
and some interesting identities for 3-cocycles appear when
naturality of the action is verified.
\end{abstract}

\maketitle

\section{Introduction}

The theory of subfactors (\cite{J1}) has experienced several new
developments in
the last few years through the introduction of
planar algebra technology (\cite{J2}).
Every subfactor comes with a very rich mathematical
object, the {\it standard invariant} or {\it planar algebra}
of the subfactor, which in nice situations is a complete
invariant of the subfactor (\cite{P2}, \cite{P3}). It can be
described in many interesting ways, as for instance a
certain category of bimodules (\cite{Oc2}, see also \cite{Bi2}),
as lattices of multi-matrix algebras (\cite{P1}),
or as a planar algebra (\cite{J2}). The planar algebra approach is
particularly
powerful since it allows one to use algebraic-combinatorial methods
in conjunction with topological ones to investigate the structure
of subfactors. A number of examples of explicit planar algebras associated to
subfactors have been computed (see for instance
\cite{BDG}, \cite{BJ}, \cite{G}, \cite{J2}, \cite{KLS}, \cite{L},
\cite{MPS}) but there is a need for more concrete examples.
This is what we accomplish in this paper. We give a description of
the planar algebra of the diagonal subfactors associated to
a $G$-kernel.

Let $P$ be a II$_1$ factor and let $\theta_1$, $\dots$, $\theta_n$
be automorphisms of $P$ (we may assume without loss of generality
that $\theta_1 = id$). Consider the subfactor
$N = \{ \sum_{i=1}^n \theta_i(x) e_{ii} \, | \, x \in P\} \subset
M = P \otimes M_n(\mathbb C)$, where $(e_{ij})_{1\le i,j \le n}$ denote
matrix units in $M_n(\mathbb C)$. $N \subset M$ is then called
the {\it diagonal
subfactor} associated to $\{\theta_i\}_{1 \le i \le n}$.
These subfactors were proposed by Jones in 1985 to provide examples
of potentially non-classifiable subfactors, since this construction
allows one to translate problems on classification of group actions
into problems on subfactors. Popa used them to prove
vanishing of 2-cohomology results for cocycle actions of finitely
generated, strongly amenable groups on an arbitrary II$_1$ factor
(\cite{P5}). Ocneanu had proved such a result for cocycle actions
of amenable groups on the hyperfinite II$_1$ factor using different
techniques (\cite{Oc}). The diagonal subfactors are of course reducible
and have Jones index $n^2$. They
provide a wealth of simple examples of infinite depth subfactors whose
structure theory is well understood. In particular, the standard
invariant or planar algebra of these subfactors has been determined
in (\cite{P2}, \cite{J2}, \cite{Bi1}).

Let $G$ be the group generated by the images $g_i$ of
$\theta_i$, $1 \le i \le n$, in $\Out P = \Aut P / \Int P$.
Popa showed that analytical properties of these subfactors
are reflected in the corresponding properties of the group $G$.
For instance, if $P$ is hyperfinite, then the diagonal
subfactor is amenable (in the sense of Popa) if and only
if $G$ is an amenable group (\cite{P2}). The subfactor has
property (T) in the sense of Popa if and only if
$G$ has property (T) of Kazhdan (\cite{P4}). The principal
graphs of these subfactors are then Cayley-like graphs of $G$
with respect to the generators $g_1$, $\dots$,
$g_n$ and their inverses (see \cite{P2} or \cite{Bi1} for the precise
statement). The higher relative commutants, Jones projections and
conditional expectations have all been worked out
in (\cite{P2}, \cite{J2}, \cite{Bi1}).

A well-known variation of the diagonal subfactors is obtained
as follows (see e.g. \cite{P2}). Consider a $G$-kernel, that is
an injective
homomorphism $\chi$ from a (countable, discrete, finitely generated)
group $G$ into $\Out P$. Denote by
$\epsilon : \Aut P \to \Out P$ the canonical homomorphism,
and let $\alpha: G \to \Aut P$ be a lift of $\chi$ such
that $\epsilon \circ \alpha_{s} = \chi(s)$, for all $s \in G$.
It follows that $\alpha_s \alpha_t = \Ad u(s,t) \alpha_{st}$,
for all $s$, $t \in G$, and some unitaries $u(s,t) \in P$. Associativity
of composition of automorphisms implies that
$\Ad( u(r,s)u(rs,t)) = \Ad (\alpha_r(u(s,t))u(r,st))$. Hence there
is a function $\omega : G \times G \times G \to \mathbb T$ with
$u(r,s)u(rs,t)=\omega(r,s,t)\alpha_r(u(s,t))u(r,st)$. It is easy
to check that $\omega$ is a 3-cocycle and that its class in
$H^3(G, \mathbb T)$ does not depend on the choices made. One usually
denotes the class by $\Ob(G)$ or $\Ob(\chi)$, the {\it obstruction}
of $\chi$. It is an obstruction to
lifting $G$ to an {\it action} on the II$_1$ factor $P$. Clearly, if
two $G$-kernels $\chi$ and $\eta$ are conjugate (in $\Out P$), then
$\Ob(\chi)=\Ob(\eta)$. It was shown in \cite{Oc} that $\Ob$
is a complete conjugacy invariant for $G$-kernels if $P$ is the
hyperfinite II$_1$ factor and $G$ is a countable, discrete, amenable
group. Note that in general, even if $\Ob(G)=0$, there may be no lifting
of the $G$-kernel to an action on $P$. Connes and Jones found in
\cite{CJ} the first such example of a non-liftable $G$-kernel
with vanishing obstruction, where $G$ is a group with property (T).
Vanishing of the obstruction is a necessary and sufficient
condition for $G$ to lift to a {\it cocycle action} on the
II$_1$ factor $P$. See \cite{BNP}, \cite{J3}, \cite{Oc}, \cite{P5}, \cite{Su}
for more on this.

Connes showed that if $G$ is cyclic, one can construct $G$-kernels
in $\Out R$ with arbitrary obstructions, where $R$ denotes the
{\it hyperfinite} II$_1$ factor (\cite{Co2}). It was an open problem
whether this result is true for more general groups, and Jones
showed in \cite{J0} that this is indeed the case for $G$
an arbitrary discrete group. Thus,
given a discrete group $G$ and a class $\pi \in H^3(G, \mathbb T)$,
there is $G$-kernel $\chi : G \to \Out R$ with $\Ob(\chi) = \pi$.
Sutherland constructed $G$-kernels with arbitrary obstructions in
non-hyperfinite II$_1$ factors (\cite{Su}).

Given a finitely generated, countable group $G$ and a 3-cocycle $\omega$,
we can associate a (hyperfinite) subfactor to $(G, \omega)$ as
follows. Let
$\chi : G \to \Out P$ be a $G$-kernel with $\Ob(G) = [\omega] \in
H^3(G, \mathbb T)$ (choose for instance $P=R$, and use \cite{J0}).
Fix generators $\{ g_1, \dots , g_n \}$
of G, let $\alpha$ be a lift of $\chi$,  and consider the diagonal
subfactor associated to the
automorphisms $\alpha_{g_1}$, $\dots$, $\alpha_{g_n}$ (let us choose
$g_1=e$, the identity of $G$, and $\alpha_e = id$). If $\eta$ is
another $G$-kernel with lift $\beta$, and automorphisms
$\beta_{g_1}$, $\dots$, $\beta_{g_n}$, then the diagonal subfactors
associated to $(\alpha_{g_i})_i$ and $(\beta_{g_j})_j$ are
isomorphic if and only if there is an automorphism $\theta$ of
the underlying II$_1$ factor $P$ such that $\alpha_{g_{\pi(i)}}
= \theta \beta_{g_i} \theta^{-1} \mod \Int P$, where $\pi$ is
a permutation of the indices (fixing $1$) (see e.g. \cite{P5}).
Thus, in particular, isomorphism of these diagonal subfactors implies
that $\Ob(\chi) = \Ob(\eta)$ (up to a possible modification of
$\chi$ by the permutation $\pi$). Isomorphism of standard invariants
is weaker than isomorphism of subfactors, but we still have the
following. If $(G,\chi)$ and $(G, \eta)$ are
two $G$-kernels as above with $\Ob(\eta) = \Ob(\chi)$,
then the standard invariants of the associated diagonal subfactors
are isomorphic (see (\cite{P2}, page 228 ff.). The converse is not
true due to the fact that group isomorphisms can change the class
of a 3-cocycle. If one constructs diagonal subfactors where the
automorphisms are repeated with distinct multiplicities, one can
show a converse, e.g. if $G$ is strongly amenable and $P$ is hyperfinite,
using Popa's classification results of amenable subfactors
(\cite{P2}, page 229 ff.). Clearly, the 3-cocycle $\omega$
giving rise to the obstruction of the $G$-kernel will appear in the standard
invariant of these diagonal subfactors (with cocycle) and the purpose
of this paper is to give a precise description of this occurrence.

Here is a more detailed outline of the sections of this paper.
We review group cocycles in section 2. In
section 3, we define an abstract planar algebra
$P^{\langle g_i : i \in I \rangle, \omega}$
associated to a finitely generated group $G$ with a fixed
finite generating set $\{ g_i \}_{i \in I}$, and a 3-cocycle
$\omega \in Z^3(G,\mathbb T)$. The vector spaces underlying this
planar algebra are spanned by multi-indices in $I^{2n}$
such that the corresponding alternating word on generators and
their inverses is the identity in $G$. The action of Jones' planar
operad is defined explicity, and of course $\omega$ appears
prominently in this definition. It should be noted that the
definition of the action of a tangle involves a labelling of
the {\it strings} in the tangle whereas the planar algebra
description of the group-type subfactors in \cite{BDG} involved
a labelling of {\it boundary segments} (called ``openings'' in
\cite{BDG}) of the internal and external discs of the tangle.
We would like to point out that the 3-cocycle $\omega$ does not
appear in our definition of
the action of the multiplication, inclusion, Jones projection
and right conditional expectation tangles. It does appear in
the definition of the action of the left conditional expectation
tangles (and hence the rotation tangles). We verify in this section
that our definition of the action of tangles is indeed natural
with respect to composition of tangles. This takes a little work,
but some interesting identities for 3-cocycles appear along the way.

In section 4 we give a model for the higher relative commutants
of the diagonal subfactor with cocycle and describe the associated
concrete planar algebra. We choose an appropriate basis of the
higher relative commutants which allows us to identify this
concrete planar algebra with the abstractly defined
one in section 3. This isomorphism is obtained in the usual way
by constructing a filtered $*$-algebra isomorphism between
the abstract planar algebra of section 3 and the concrete one
of section 4, that preserves Jones projection and
conditional expectation tangles. The main feature of this
planar algebra is the fact that the distinguished basis
of the higher relative commutants that we found here, matches
with the one coming from the description of the
planar algebra as a path algebra associated to the principal
graphs of the subfactor (see e.g. \cite{JS}, \cite{J2}).
Conversely, we prove that any finite index extremal 
subfactor whose standard invariant is given by the 
abstract planar algebra (in Section 3), must necessarily 
be (isomorphic to) a diagonal subfactor. 

\section{A brief review of group cocycles}
For the convenience of the reader, we recall in this brief section
the definition of a cocycle of a group $G$. $G$ will denote
throughout this article a countable, discrete group, and we
will denote the identity of $G$ by $e$. Define $C^n =
Fun(G^n,\mathbb T)$ the space of functions from $G^n$ to $\mathbb T$ and
$\partial^n:C^n \rightarrow C^{n+1}$ by\\
\begin{align*}
& \partial^n(\phi) (g_1, \cdots, g_{n+1})\\
= & \; \phi(g_2, \cdots, g_{n+1}) \; {\phi(g_1 g_2, g_3, \cdots,
g_{n+1})}^{-1} \; \phi(g_1, g_2 g_3, g_4, \cdots, g_{n+1}) \cdots\\
& \; \cdots {\phi(g_1, \cdots, g_{n-1}, g_n g_{n+1})}^{(-1)^n} \;
{\phi(g_1, \cdots, g_n)}^{(-1)^{n+1}}
\end{align*}
It follows that $\left(\partial^{n+1} \circ \partial^n\right) (\cdot)
= 1_{C^{n+2}}$ where $1$ denotes the constant function $1$. Denote
$ker(\partial)$ by $Z^n(G,\mathbb T)$ (whose elements are called
$n$-cocycles) and $Im(\partial^{n-1})$ by $B^n(G,\mathbb T)$ (whose elements
are called coboundaries). Note that
$B^n(G,\mathbb T) \subset Z^n(G,\mathbb T)$.

In this paper we will be dealing mostly with a $3$-cocycle $\omega$.
Thus $\omega$ satisfies the identity

\begin{align}\label{defcocycle}
\omega(g_1,g_2,g_3)\;\omega(g_1,g_2g_3,g_4)\;\omega(g_2,g_3,g_4) =
\omega (g_1g_2,g_3,g_4)\;\omega(g_1,g_2,g_3g_4)
\end{align}

We call $\omega$ {\it normalized} if $\omega(g_1,
g_2, g_3) = 1$ whenever either of $g_1$, $g_2$, $g_3$ is $e$.

Any cocycle $\omega$ is coboundary equivalent to a normalized cocycle.
In particular, $(\omega \cdot \partial^2(\phi))$ is a normalized
$3$-cocycle where $\phi \in C^2$ is defined as $\phi(g_1,g_2) =
\omega(g_1,e,e) \; \overline{\omega}(e,e,g_2)$ for all $g_1, g_2 \in G$.
\section{An abstract planar algebra}\label{plnalg}
In this section we will construct an abstract planar algebra
which, in section \ref{subfact}, will be shown to be isomorphic
to the planar algebra of a diagonal subfactor with cocycle.

Let $G$ be a countable, discrete group
generated by a finite subset $\{ g_i \}_{i \in I}$, and let
$\omega \in Z^3(G,\mathbb T)$ be a normalized $3$-cocycle.
We will construct a planar algebra
$P^{\langle g_i : i \in I \rangle, \omega}$ (= $P$) as follows.
Let $e$ denote the identity of $G$. We define first a
map $alt$ from multi-indices $\coprod_{n \ge 0}I^n$ to $G$ by
\[ \left( \coprod_{n \geq 0} I^n \right) \ni \ul{i} = (i_1, \cdots ,
  i_n) \stackrel{alt}{\longmapsto} g^{-1}_{i_1} g_{i_2} \cdots
  g^{(-1)^n}_{i_n} =   alt(\ul{i}) \in G \]
where $alt$ of the empty multi-index is defined to be the identity
element of the group.
To define the planar algebra we need to define vector spaces
$P_n$ and an action of Jones' planar operad on these vector spaces.
We refer to \cite{J2} for the planar algebra terminology used in
this paper.

\subsection*{The vector spaces} For $n \geq 0$, define
$P_n = \left
  \{
\begin{array}{ll}
\C \{\ul{i} \in I^{2n} : alt(\ul{i}) = e \} & \text{if $n > 0$},\\
\C & \text{if $n = 0$.}
\end{array}
\right. $

\subsection*{Action of tangles}
\begin{figure}
\begin{center}
\epsfig{file=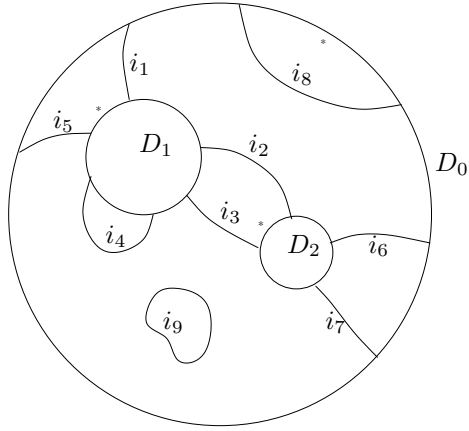, scale=0.4}
\caption{Example of a state on a tangle}
\label{fig:diagstate}
\end{center}
\end{figure}

Let $T$ be an $n_0$-tangle having
internal discs $D_1, \cdots,  D_b$ with colors $n_1, \cdots, n_b$
respectively (or no internal discs of course). A {\it state}
$\sigma$ on $T$ is a map from
$\{ \text{strings in $T$} \}$ to $I$ such that $alt({\sigma
  |_{\partial D_c}}) = e$ for all $1 \leq c \leq b$ where $\sigma
|_{\partial D_c}$ denotes the element of $I^{2 n_c}$ obtained by reading
  the elements of $I$ induced at the marked points on the
  boundary of $D_c$ by the strings via the map $\sigma$, starting
  from the first marked point and moving clockwise. This has been
  illustrated in Figure \ref{fig:diagstate}, where $alt({\sigma
  |_{\partial D_1}})$ and $alt({\sigma|_{\partial D_2}})$ are just the
 products  $g^{-1}_{i_1}g_{i_2}g^{-1}_{i_3}g_{i_4}g^{-1}_{i_4}g_{i_5}$
and $g^{-1}_{i_2}g_{i_6}g^{-1}_{i_7}g_{i_3}$ respectively, and are
thus required to be the identity element. It is a consequence that
  $alt({\sigma |_{\partial D_0}}) (=
g^{-1}_{i_8}g_{i_6}g^{-1}_{i_7}g_{i_5}g^{-1}_{i_1}g_{i_8} = e$
in the figure) holds for the external disc.
Let $\mathcal{S} (T)$ denote the states on $T$.

\psfrag{D0}{$D_0$}
\psfrag{D1}{$D_1$}
\psfrag{D2}{$D_2$}
\psfrag{i1}{$i_1$}
\psfrag{i2}{$i_2$}
\psfrag{i3}{$i_3$}
\psfrag{i4}{$i_4$}
\psfrag{i5}{$i_5$}
\psfrag{i6}{$i_6$}
\psfrag{i7}{$i_7$}
\psfrag{i8}{$i_8$}
\psfrag{i9}{$i_9$}

In order to define the action $Z_T$ of $T$, it is enough to define the
coefficient $\langle Z_T (\ul{k}^1, \cdots, \ul{k}^b) \: | \: \ul{k}^0
\rangle$ of $\ul{k}^0$ in the linear expansion of $Z_T (\ul{k}^1,
\cdots, \ul{k}^b)$ where $\ul{k}^c \in I^{2n_c}$ such that
$alt({\ul{k}^c}) = e$ for $1 \leq c \leq b$. For this, we choose a
picture $T_1$ in the isotopy class of $T$ and then choose a simple
path $p_c$ in $D_0 \setminus [\coprod^b_{c=1} Int(D_c)]$ starting from
the $*$
\footnote{Recall $*$ of a disc $D$ is a point chosen on the boundary
of $D$ strictly between the last and the first marked points, moving
clockwise.} of $D_0$ to that of $D_c$ for $1 \leq c \leq b$ such that:

(i) $p_c$ intersects the strings of $T_1$ transversally for $1 \leq c
\leq b$,

(ii) $p_{c_1}$ and $p_{c_2}$ intersects exactly at the $*$ of $D_0$
for $1 \leq c_1 \neq c_2 \leq b$.\\
Note that any state $\sigma$ on $T$ gives an element $\sigma |_{p_c}
\in I^{m_c}$ obtained by reading the elements of $I$ induced by
$\sigma$ at the crossings of the path $p_c$ and the strings along the
direction of the path where $m_c$ (necessarily even) is the number of
strings cut by $p_c$.

Define
\[\langle Z_T (\ul{k}^1, \cdots, \ul{k}^b) \: | \: \ul{k}^0
\rangle \; = \sum_{\substack{\sigma \in \mathcal{S}(T) \text{ s.t.}\\
\sigma |_{\partial D_d} = \ul{k}^d \\ \text{for } 0 \leq d \leq b}} \;
\prod^b_{c=1} \; \lambda_{\sigma|_{p_c}} (\ul{k}^c)\]
where $\lambda_{\ul{j}} (\ul{i}) = \prod^n_{s=1} \lambda_{\ul{j}}
(\ul{i},s)$ and

$$\lambda_{\ul{j}} (\ul{i},s) = \left\{
\begin{array}{l}
\ol{\omega} (alt({\ul{j}}) , alt(i_1, \cdots, i_s) ,
g_{i_s}) \text{ if $s$ is odd}\\
\omega (alt({\ul{j}}) , alt(i_1, \cdots, i_{s-1}) , g_{i_s}) \text{
  if $s$ is even} \end{array} \right.$$

for $\ul{i} \in I^n$, $\ul{j} \in I^m$. If there is no
compatible state on $T$, then we take the coefficient to be $0$ and if
there is no internal disc in $T$, then the scalar inside the sum is
considered to be $1$. (Note that
$\lambda_{\ul{j}} (\ul{i})$ depends only on $alt({\ul{j}})$ and
$\ul{i}$.)

We need to show first that the multi-linear map $Z_T$ is
well-defined. Two configurations of paths $\{p_c\}^b_{c=1}$ and
$\{p^{\prime}_c\}^b_{c=1}$ in $T_1$ can be obtained from
each other using a finite sequence of the following moves:\\
\psfrag{D0}{$D_0$}
\psfrag{D1}{$D_1$}
\psfrag{D2}{$D_2$}
\psfrag{Di}{$D_c$}
\psfrag{Dj}{$D_d$}
\psfrag{p}{$p$}
\psfrag{p1}{$p_1$}
\psfrag{p2}{$p_2$}
\psfrag{pi}{$p_c$}
\psfrag{pj}{$p_d$}
\psfrag{c}{$\cdots$}
\psfrag{A}{$A$}
\psfrag{B}{$B$}
\psfrag{*}{$*$}
\psfrag{equiv}{$\sim$}
\psfrag{1to2r}{$1~\cdots~2r$}
\psfrag{1}{$1$}
\psfrag{m1}{$m_1$}
\psfrag{m1+1}{$m_1+1$}
\psfrag{2m1}{$2m_1$}
\psfrag{n1}{$n_1$}
\psfrag{n1+1}{$n_1+1$}
\psfrag{2n1}{$2n_1$}
\psfrag{n2}{$n_2$}
\psfrag{n2+1}{$n_2+1$}
\psfrag{2n2}{$2n_2$}
I. {\it isotopy}
\begin{center}
      \epsfig{file=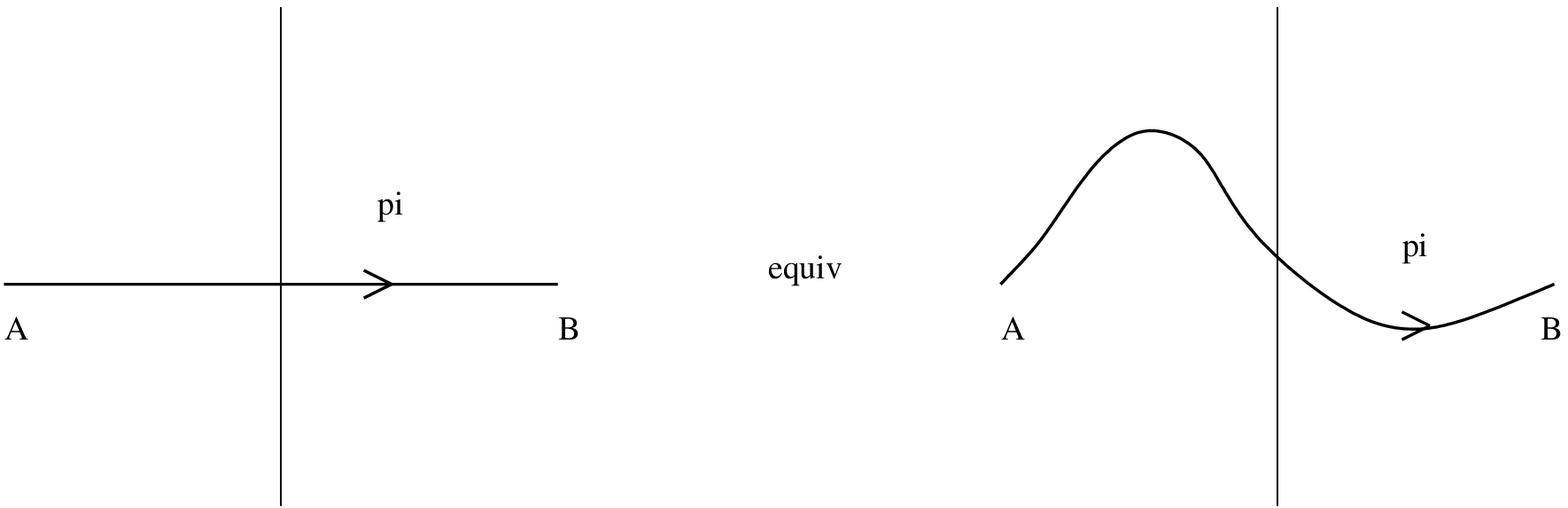, scale=0.5}
    \end{center}
II. {\it cap-sliding moves}
\begin{center}
      \epsfig{file=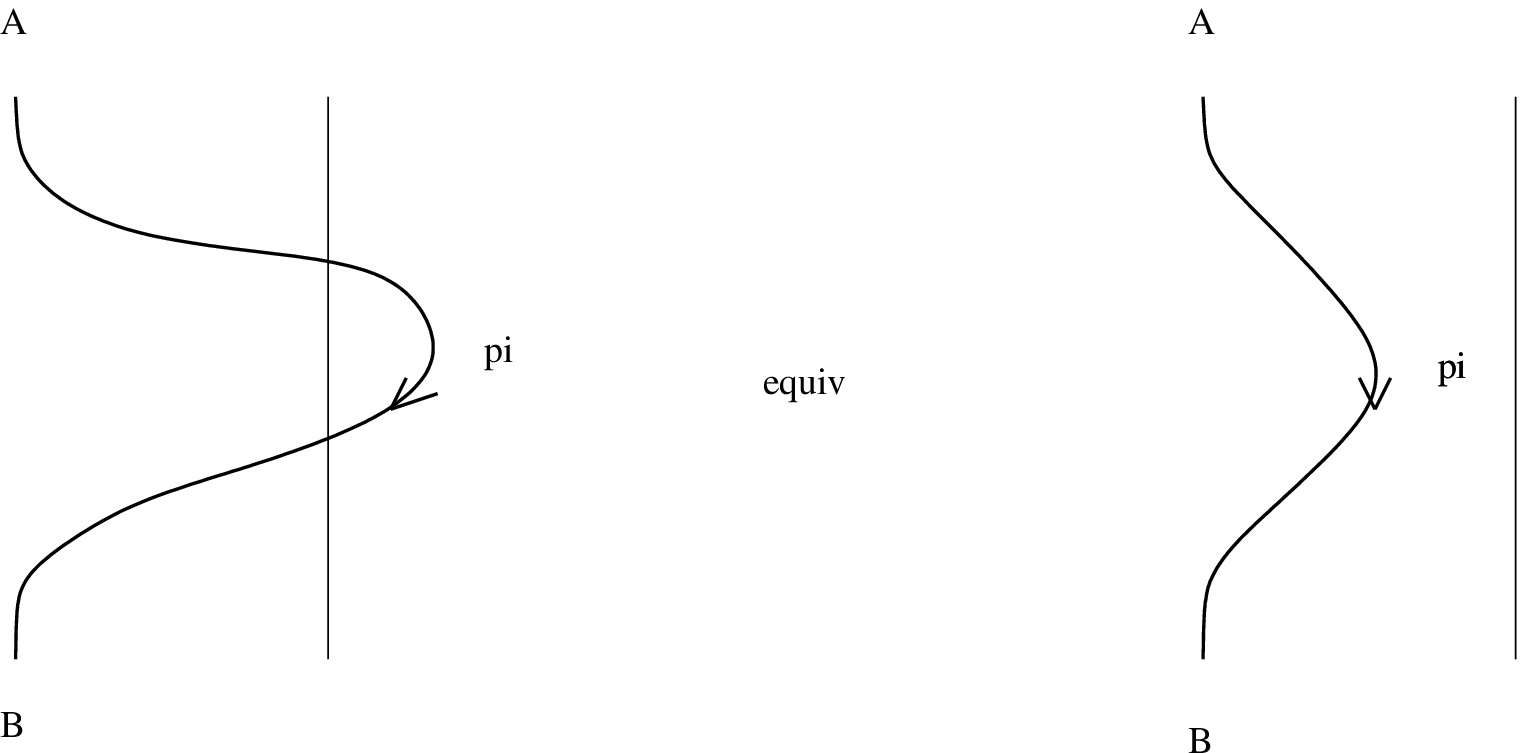, scale=0.5}
    \end{center}
III. {\it disc-sliding moves}
\begin{flushleft}
      \epsfig{file=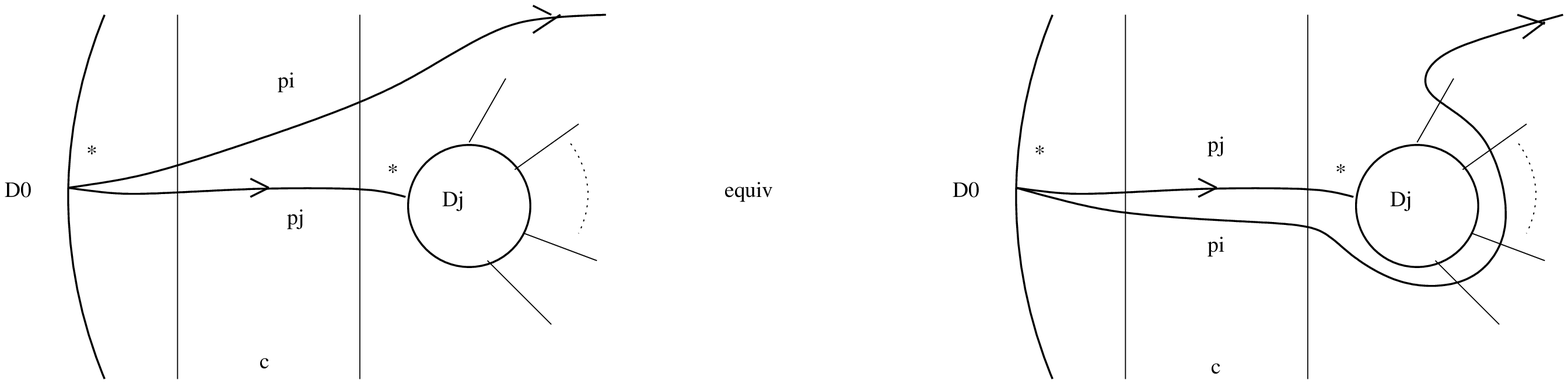, scale=0.5}
    \end{flushleft}
IV. {\it rotation moves}
\begin{flushleft}
      \epsfig{file=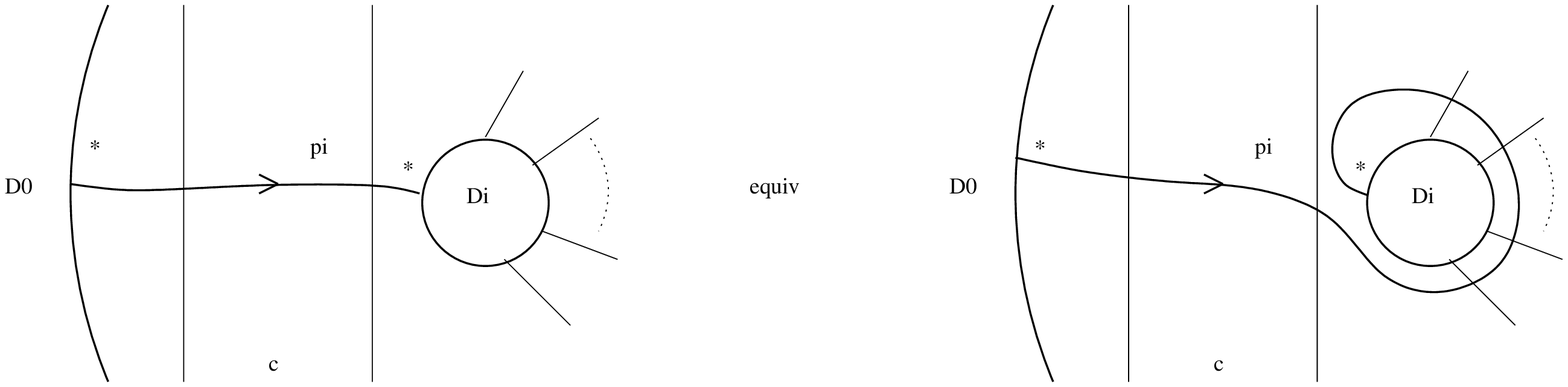, scale=0.5}
    \end{flushleft}
It is enough to check that the definition of the action is independent
under each of the
above moves. Invariance under isotopy moves are the easiest to check
since $\sigma |_{p_c} = \sigma |_{p^{\prime}_c}$ for all $c \in \{1,
\cdots, b\}$. To see invariance under the remaining three moves, we show that
$alt({\sigma |_{p_c}}) = alt({\sigma |_{p^{\prime}_c}})$ for
$1 \leq c \leq b$. For a cap-sliding move, note that the cap induces
the same index at
the two consecutive crossing with the path but after applying the $alt$
map, the corresponding group elements cancel each other since they
inverses of each other. For the disc-sliding (resp., rotation moves),
the invariance follows from the fact that
$alt({\sigma |_{\partial{D_d}}}) = e$ (resp.,
$alt({\sigma |_{\partial{D_c}}}) = e $).

\medskip
We will show next that the action is compatible with composition
of tangles.

\subsection*{Action is natural with respect to composition of tangles.}
Let $S$ be an $m_0$-tangle
containing the internal discs $D^{\prime}_1, \cdots, D^{\prime}_a$
with colors $m_1, \cdots, m_a$ and $T$ be an $m_1$-tangle containing
internal discs $D_1, \cdots, D_b$ with colors $n_1, \cdots, n_b$. Let
$D^{\prime}_0$ and $D_0$ denote the external discs of $S$ and $T$
respectively. We
need to show that for all $\ul{i}^{c^\prime} \in I^{2m_{c^{\prime}}}$
and $\ul{j}^c \in I^{2n_c}$ where $c^\prime \in \{0, 2, \cdots, a\}$
and $c \in \{1, \cdots, b\}$
\begin{equation}\label{comp}
\langle Z_S ( Z_T (\ul{j}^1, \cdots, \ul{j}^b ), \ul{i}^2,
\cdots, \ul{i}^a) | \ul{i}^0 \rangle =
\langle Z_{S
\circ_{D^{\prime}_1} T} (\ul{j}^1, \cdots, \ul{j}^b, \ul{i}^2,
\cdots, \ul{i}^a) | \ul{i}^0 \rangle
\end{equation}
The left-hand side of (\ref{comp}) can be expanded as
\begin{eqnarray}
&& \sum_{\substack{\ul{i}^1 \in I^{2m_1} \text{ s.t.}\\
alt({\ul{i}^1}) = e}} \langle Z_S ( \ul{i}^1,
\ul{i}^2, \cdots, \ul{i}^a) | \ul{i}^0 \rangle \; \langle Z_T (\ul{j}^1,
\cdots, \ul{j}^b ) | \ul{i}^1 \rangle \nonumber\\
&=& \sum_{\substack{\sigma \in \mathcal{S}(S), \: \tau \in
\mathcal{S}(T) \text{ s.t. } \sigma |_{\partial{D^{\prime}_{c^\prime}}} =
\ul{i}^{c^\prime} , \: \tau|_{\partial{D_c}} = \ul{j}^c \text{ for} \\
c^\prime \in \{0,2, \cdots, a\}, \: c
 \in \{1, \cdots, b\} \text{ and } \: \sigma |_{\partial{D^{\prime}_1}}
= \tau|_{\partial{D_0}} = \ul{i}^1 }} \left( \prod^a_{c^{\prime} =1} \lambda_{
\sigma
 |_{p^{\prime}_{c^\prime}}} (\ul{i}^{c^\prime})\right) \left(
 \prod^b_{c=1} \lambda_{\tau |_{p_c}} (\ul{j}^c) \right) \nonumber
\end{eqnarray}
where $p^{\prime}_{c^\prime}$'s and $p_c$'s are paths in the tangles
$S$ and $T$ respectively, required to define their actions. For the
action of $S \circ_{D^{\prime}_1} T$, we consider the paths
$p^{\prime}_{c^\prime}$ for $2 \leq c^\prime \leq a$ and $p^{\prime}_1
\circ p_c$ for $1 \leq c \leq b$. (Strictly speaking, one has to
disjointify the $p^{\prime}_1$-portion of the paths $(p^{\prime}_1
\circ p_c)$ for different values of $c$ in order to define the action
of $S \circ_{D^{\prime}_1} T$.) So, the right-hand side of (\ref{comp})
becomes
\[
\sum_{\substack{\gamma \in \mathcal{S}(S \circ_{D^{\prime}_1} T)
\text{ s.t. } \gamma |_{\partial{D^{\prime}_{c^\prime}}} =
\ul{i}^{c^\prime} , \: \gamma
|_{\partial{D_c}} = \ul{j}^c\\ \text{for } c^\prime \in \{0,2, \cdots, a\}, \:
c
 \in \{1, \cdots, b\} }} \left( \prod^a_{c^{\prime} =2} \lambda_{\gamma
 |_{p^{\prime}_{c^\prime}}} (\ul{i}^{c^\prime})\right) \left(
 \prod^b_{c=1} \lambda_{\gamma |_{p^\prime_1 \circ p_c}} (\ul{j}^c)
 \right)
\]
Observe that $
\left\{
(\sigma, \tau) \in \mathcal{S}(S) \times \mathcal{S}(T) :
\begin{array}{c}
\sigma |_{\partial{D^{\prime}_{c^\prime}}} = \ul{i}^{c^\prime} \text{ for }
c^\prime \in \{0,2, \cdots, a\} \\
\tau|_{\partial{D_c}} = \ul{j}^c \text{ for } 1 \leq c \leq b,
 \: \sigma |_{\partial{D^{\prime}_1}} = \tau |_{\partial{D_0}}
\end{array}
\right\}$
is clearly in bijection with $
\left\{
\gamma \in \mathcal{S}(S \circ_{D^{\prime}_1} T) :
\begin{array}{c}
\gamma |_{\partial{D^{\prime}_{c^\prime}}} = \ul{i}^{c^\prime} \text{ for }
c^\prime \in \{0,2, \cdots, a\},\\ \gamma |_{\partial{D_c}} = \ul{j}^c \text{
 for } 1 \leq c \leq b
\end{array}
\right\}$.
A bijection is obtained by sending $(\sigma, \tau)$ to the state defined by
$\sigma$ (resp. $\tau$) on the $S$-part (resp. $T$-part) of $S
\circ_{D^{\prime}_1} T$ and the well-definedness of such a state is a
consequence of the condition $\sigma |_{\partial{D^{\prime}_1}} =
\tau |_{\partial{D_0}}$;
we denote this state by $\sigma \circ \tau$. If these sets are empty,
equation (\ref{comp}) holds trivially since both sides have value
$0$. Let us assume that the sets are nonempty. It is enough to prove
for $\sigma \in \mathcal{S}(S)$ and $\tau \in \mathcal{S}(T)$ such
that $\sigma |_{\partial{D^{\prime}_{c^\prime}}} = \ul{i}^{c^\prime}$ for $0
\leq c^\prime \leq a$ and $\tau |_{\partial{D_c}} = \ul{j}^c$ for
$1 \leq c \leq
b$ and $\tau |_{\partial{D_0}} = \ul{i}^1$ we have
\begin{equation}\label{comp2}
\lambda_{\sigma
 |_{p^{\prime}_1}} (\ul{i}^1)
 \prod^b_{c=1} \lambda_{\tau |_{p_c}} (\ul{j}^c)
\; = \;
\prod^b_{c=1} \lambda_{(\sigma \circ \tau) |_{(p^\prime_1 \circ p_c)}}
 (\ul{j}^c)
\end{equation}
We prove this in two cases.\\

\medskip
\noindent
{\bf Case 1:} $T$ has no internal disc or closed loop, that is, $T$ is
a Temperley-Lieb diagram. Then the right-hand side of equation
(\ref{comp2}) is
$1$. It remains to show $\lambda_{\sigma |_{p^{\prime}_1}} (\ul{i}^1)
= 1$. This follows from the next lemma and the fact $\ul{i}^1 = Z_T
(1)$ is a sequence of non-crossing matched pairings of indices from $I$.

\begin{lemma}\label{pair1}
If $i \in I$, $\ul{j} \in I^m$, $\ul{i} = (i_1, \cdots, i_{n}) \in
I^{n}$ and $0 \leq s \leq n$, then we have $\lambda_{\ul{j}}
(\ul{i}) = \lambda_{\ul{j}} (i_1, \cdots, i_s, i, i, i_{s+1}, \cdots,
i_{n})$.
\end{lemma}

\noindent{\bf Proof:} Note that\\
(i) $\lambda_{\ul{j}} (\ul{i} , r) = \lambda_{\ul{j}} ((i_1, \cdots,
i_s, i, i, i_{s+1}, \cdots, i_{n}),r)$ for $1 \leq r \leq s$, and\\
(ii)  $\lambda_{\ul{j}} (\ul{i} , r) = \lambda_{\ul{j}} ((i_1, \cdots,
i_s, i, i, i_{s+1}, \cdots, i_{n}),r+2)$ for $s+1 \leq r \leq n$.\\

We compute then,
\begin{align*}
& \lambda_{\ul{j}} ((i_1, \cdots, i_s, i, i, i_{s+1}, \cdots,
  i_n),s+1) \; \lambda_{\ul{j}} ((i_1, \cdots, i_s, i, i, i_{s+1},
  \cdots, i_n),s+2) \\
= &
\left\{\begin{array}{ll}
\ol{\omega}
(alt({\ul{j}}) , alt(i_1, \cdots, i_s, i), g_i) \; \omega
  (alt({\ul{j}}) , alt(i_1, \cdots, i_s, i) , g_i) & \text{if $s$ is
  even}\\
\omega (alt({\ul{j}}) , alt(i_1, \cdots, i_s) , g_i) \; \ol{\omega}
  (alt({\ul{j}}) , alt(i_1, \cdots, i_s, i, i) , g_i)  & \text{if
  $s$ is odd}
\end{array}\right.\\
= & 1.
\end{align*}\qed

\begin{remark}\label{pair2}
In Lemma \ref{pair1}, if $\ul{i}$ is a sequence of indices with
non-crossing matched pairings, then we can apply the lemma
several times to reduce all the consecutive matched pairings
to get $\lambda_{\ul{j}} (\ul{i}) = 1$.\qed
\end{remark}

\medskip
\noindent
{\bf Case 2:} $T$ has at least one internal disc. Any unlabelled
tangle $T$ can be expressed as composition of elementary annular
tangles of four types as described in \cite{BDG}, namely, capping,
cap-inclusion, left-inclusion and disc-inclusion tangles. It is enough
to prove equation (\ref{comp2}) for any tangle $S$ and compatible
tangle $T$ in
$\mathcal{E}$ ($=$ the set of all elementary tangles). If $T$ is of
capping or cap-inclusion type annular tangle, the proof directly
follows from Lemma \ref{pair1} and is left to the reader.

If $T$ is a left-inclusion annular tangle, equation (\ref{comp2}) is
implied from the following lemma.

\begin{lemma}\label{leftinc1}
For all $\ul{i} \in I^{2m}$, $\ul{j} \in I^{2n}$ and $\ul{k} \in I^{2n_1}$
such that $alt({\ul{i}}) = e$, we have
\[
\lambda_{(\ul{k}, \ul{j})} (\ul{i}) =
\lambda_{\ul{k}} (\ul{j}, \ul{i}, \widetilde{\ul{j}}) \; \lambda_{\ul{j}}
(\ul{i})
\]
where $\widetilde{\ul{j}}$ is the sequence of indices from $\ul{j}$ in the
reverse order.
\end{lemma}

\noindent{\bf Proof:} We rearrange the terms of the right-hand side of the
identity in the lemma in the following way:
\begin{align*}
& \lambda_{\ul{k}} (\ul{j}, \ul{i}, \widetilde{\ul{j}}) \;
  \lambda_{\ul{j}} (\ul{i})\\
= & \left(\prod^{2n}_{r=1} \lambda_{\ul{k}} ((\ul{j}, \ul{i},
\widetilde{\ul{j}}),r) \; \lambda_{\ul{k}} ((\ul{j}, \ul{i},
\widetilde{\ul{j}}), 2(m+2n)-r+1)\right)
\left(
\prod^{2m}_{s=1} \lambda_{\ul{k}} ((\ul{j}, \ul{i},
  \widetilde{\ul{j}}) , 2n+s) \; \lambda_{\ul{j}} (\ul{i} , s)
\right)
\end{align*}
Let $\ul{i} = (i_1, \cdots, i_{2m})$ and $\ul{j} = (j_1, \cdots,
j_{2n})$. Note that for $1 \leq r \leq 2n$,\\
\begin{align*}
& \lambda_{\ul{k}} ((\ul{j}, \ul{i}, \widetilde{\ul{j}}),
  2(m+2n)-r+1)\\
= & \left\{\begin{array}{ll}
\omega(alt({\ul{k}}) , alt(j_1, \cdots, j_r) ,
g_{j_r}) & \text{if $r$ is odd}\\
\ol{\omega}(alt({\ul{k}}) , alt(j_1, \cdots, j_{r-1}) , g_{j_r}) &
\text{if $r$ is even}
\end{array}
\right.\\
= & \ol{\lambda}_{\ul{k}} ((\ul{j}, \ul{i},
\widetilde{\ul{j}}),r)
\end{align*}

since $alt({\ul{i}}) = e$ and $alt(\ul{j}, \ul{i},
  \widetilde{\ul{j}}) = alt(\ul{j}) alt(\ul{i}) (alt(\ul{j}))^{-1} = e$.
Thus the first product of the terms in the rearrangement vanishes. For
the second product, if $s \in \{1, \cdots, 2m\}$ is odd, then
\begin{align*}
& \; \lambda_{\ul{k}} ((\ul{j}, \ul{i},
  \widetilde{\ul{j}}) , 2n+s) \; \lambda_{\ul{j}} (\ul{i} , s)\\
= & \; \ol{\omega}(alt({\ul{k}}) , alt({\ul{j}}) alt(i_1, \cdots, i_s)
  , g_{i_s}) \; \ol{\omega}(alt({\ul{j}}) , alt(i_1, \cdots, i_s) ,
  g_{i_s})\\
= & \; \ol{\omega}(alt({\ul{k}}) , alt({\ul{j}}) , alt(i_1, \cdots,
  i_{s-1})) \; \ol{\omega}(alt(\ul{k},\ul{j}) , alt(i_1, \cdots, i_s)
  , g_{i_s}) \cdot\\
& \; \omega(alt({\ul{k}}) , alt({\ul{j}}) , alt(i_1, \cdots, i_s))\\
= & \; \ol{\omega}(alt({\ul{k}}) , alt({\ul{j}}) , alt(i_1, \cdots,
  i_{s-1})) \; \lambda_{(\ul{k},\ul{j})} (\ul{i} , s)
  \; \omega(alt({\ul{k}}) , alt({\ul{j}}) , alt(i_1, \cdots, i_s) )
\end{align*}
where we used the defining equation (\ref{defcocycle}) of the
$3$-cocycle $\omega$ for the second equality.

Similarly, if $s \in \{1, \cdots, 2m\}$ is even, then
\begin{align*}
& \; \lambda_{\ul{k}} ((\ul{j}, \ul{i},
  \widetilde{\ul{j}}) , 2n+s) \; \lambda_{\ul{j}} (\ul{i} , s)\\
= & \; \ol{\omega}(alt({\ul{k}}) , alt({\ul{j}}) , alt(i_1, \cdots,
  i_{s-1})) \; \lambda_{(\ul{k},\ul{j})} (\ul{i} , s) \;
  \omega(alt({\ul{k}}) , alt({\ul{j}}) , alt(i_1, \cdots, i_s))
\end{align*}
Thus,
\begin{align*}
& \prod^{2m}_{s=1} \lambda_{\ul{k}} ((\ul{j}, \ul{i},
  \widetilde{\ul{j}}) , 2n+s) \; \lambda_{\ul{j}} (\ul{i} , s)
\\
=&\prod^{m}_{t=1} \left(\lambda_{\ul{k}} ((\ul{j}, \ul{i},
  \widetilde{\ul{j}}) , 2n+2t-1) \; \lambda_{\ul{j}} (\ul{i} , 2t-1) \right)
\left(  \lambda_{\ul{k}} ((\ul{j}, \ul{i},
  \widetilde{\ul{j}}) , 2n+2t) \; \lambda_{\ul{j}} (\ul{i} , 2t) \right)
\\
=& \prod^{m}_{t=1}
\left(\begin{array}{l}
\ol{\omega}(alt({\ul{k}}) , alt({\ul{j}}) ,
  alt(i_1, \cdots, i_{2t-2})) \; \lambda_{(\ul{k},\ul{j})} (\ul{i} ,
  2t-1)\\
\lambda_{(\ul{k},\ul{j})} (\ul{i} , 2t) \; \omega(alt({\ul{k}}),
  alt({\ul{j}}), alt(i_1, \cdots, i_{2t}))
\end{array}\right)\\
=& \; \ol{\omega}(alt({\ul{k}}) , alt({\ul{j}}) , e)
\left( \prod^{m}_{t=1} \lambda_{(\ul{k},\ul{j})} (\ul{i} , 2t-1) \;
\lambda_{(\ul{k},\ul{j})} (\ul{i} , 2t) \right)
  \omega(alt({\ul{k}}), alt({\ul{j}}) , alt({\ul{i}}) )\\
=& \; \lambda_{(\ul{k}, \ul{j})} (\ul{i})
\end{align*}

since $alt({\ul{i}}) = e$ and $\omega$ is normalized.\qed

\begin{figure}[h]
\begin{center}
\epsfig{file=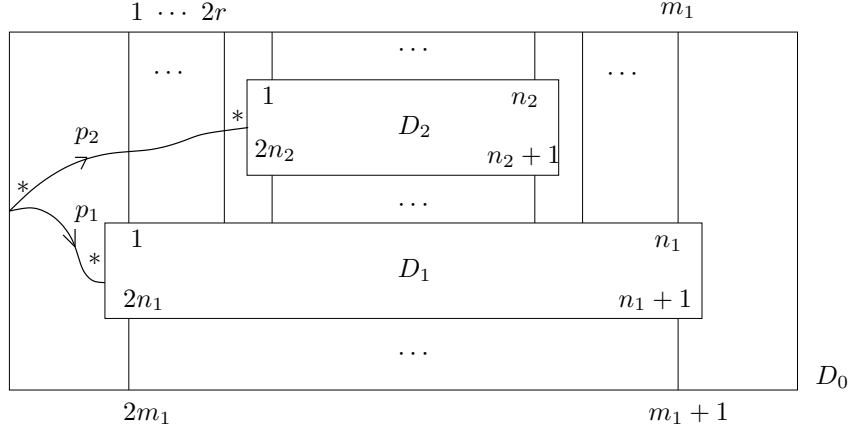, scale=0.5}
\caption{Disc inclusion tangle}
\label{path}
\end{center}
\end{figure}

\begin{remark}\label{leftinc2}
The proof of Lemma \ref{leftinc1} also implies the following identity:
\begin{align*}
\lambda_{\ul{k}} (\ul{j}, \ul{i}, \widetilde{\ul{j}})
= \prod^{2m+2n}_{s=2n+1} \lambda_{\ul{k}} ((\ul{j}, \ul{i},
  \widetilde{\ul{j}}) , s)
\end{align*}\qed
\end{remark}
Now, suppose $T$ is a disc-inclusion tangle as shown in Figure
\ref{path}. Note that $n_1 = m_1$.
Without loss of generality, we can assume that $T$ be given by the
following picture in which we also indicate the paths $p_1$ and $p_2$.

\psfrag{D0}{$D_0$}
\psfrag{D1}{$D_1$}
\psfrag{D2}{$D_2$}
\psfrag{Di}{$D_c$}
\psfrag{Dj}{$D_d$}
\psfrag{p}{$p$}
\psfrag{p1}{$p_1$}
\psfrag{p2}{$p_2$}
\psfrag{pi}{$p_c$}
\psfrag{pj}{$p_d$}
\psfrag{c}{$\cdots$}
\psfrag{A}{$A$}
\psfrag{B}{$B$}
\psfrag{*}{$*$}
\psfrag{equiv}{$\sim$}
\psfrag{1to2r}{$1~\cdots~2r$}
\psfrag{1}{$1$}
\psfrag{m1}{$m_1$}
\psfrag{m1+1}{$m_1+1$}
\psfrag{2m1}{$2m_1$}
\psfrag{n1}{$n_1$}
\psfrag{n1+1}{$n_1+1$}
\psfrag{2n1}{$2n_1$}
\psfrag{n2}{$n_2$}
\psfrag{n2+1}{$n_2+1$}
\psfrag{2n2}{$2n_2$}

Observe that since $p_1$ does not intersect any string, we have
that $\tau |_{p_1}$ is the empty multi-index. So it is enough to prove

\begin{equation}\label{discinc}
\lambda_{\sigma |_{p^{\prime}_1}} (\ul{i}^1) \;
\lambda_{\tau |_{p_2}} (\ul{j}^2)
\; = \;
\lambda_{\sigma |_{p^\prime_1}} (\ul{j}^1)
\; \lambda_{(\sigma |_{p^\prime_1}, \tau |_{p_2})} (\ul{j}^2)
\end{equation}

Let us denote $\tau |_{p_2}$ by $\ul{k} \in I^{2r}$. Since $\tau$ is a
state, the following relations clearly follow from Figure \ref{path}:\\

\noindent
(i) $i^1_s = j^1_s$ for $1 \leq s \leq 2r$ and $2r + n_2 +1 \leq s
\leq 2n_1$,\\
(ii) $i^1_{2r + s} = j^2_s$ for $1 \leq s \leq n_2$,\\
(iii) $i^1_s = k_s = j^1_s$ for $1 \leq s \leq 2r$,\\
(iv) $j^1_{2r+s} = j^2_{2n_2 - s+1}$ for $1 \leq s \leq n_2$.\\

We now express $\lambda_{\sigma |_{p^{\prime}_1}} (\ul{i}^1)$ as a
product of three terms with which we work separately.
\[
\lambda_{\sigma |_{p^{\prime}_1}} (\ul{i}^1)  =  \left(
  \prod^{2r}_{s=1} \lambda_{\sigma |_{p^{\prime}_1}}
  (\ul{i}^1, s) \right)
\left( \prod^{2r + n_2}_{s=2r+1} \lambda_{\sigma |_{p^{\prime}_1}}
  (\ul{i}^1, s) \right)
\left( \prod^{2n_1}_{s= 2r + n_2 + 1} \lambda_{\sigma |_{p^{\prime}_1}}
  (\ul{i}^1, s) \right)
\]

\noindent
First term: For $1 \leq s \leq 2r$,
\begin{align*}
\lambda_{\sigma |_{p^{\prime}_1}} (\ul{i}^1, s) & =
\left\{\begin{array}{ll}
\ol{\omega} (alt({\sigma |_{p^{\prime}_1}}), alt({i^1_1} , \cdots,
  {i^1_s}),g_{i^1_s}) & \text{if $s$ is odd}\\
\omega (alt({\sigma |_{p^{\prime}_1}}), alt({i^1_1} , \cdots,
  {i^1_{s-1}}),g_{i^1_s}) & \text{if $s$ is even}
\end{array}\right.\\
& =
\left\{\begin{array}{ll}
\ol{\omega} (alt({\sigma |_{p^{\prime}_1}}), alt({j^1_1} , \cdots,
  {j^1_s}),g_{j^1_s}) & \text{if $s$ is odd}\\
\omega (alt({\sigma |_{p^{\prime}_1}}), alt({j^1_1} , \cdots,
  {j^1_{s-1}}),g_{j^1_s}) & \text{if $s$ is even}
\end{array}\right.
\text{(applying (i))}\\
& = \; \lambda_{\sigma |_{p^{\prime}_1}} (\ul{j}^1,s).\\
\end{align*}

\noindent
Second term: For $2r+1 \leq s \leq 2r+n_2$,
\begin{align*}
\lambda_{\sigma |_{p^{\prime}_1}} (\ul{i}^1, s) & =
\left\{\begin{array}{ll}
\ol{\omega} (alt({\sigma |_{p^{\prime}_1}}), alt({i^1_1} , \cdots,
  {i^1_s}),g_{i^1_s}) & \text{if $s$ is odd}\\
\omega (alt({\sigma |_{p^{\prime}_1}}), alt({i^1_1} , \cdots,
  {i^1_{s-1}}),g_{i^1_s}) & \text{if $s$ is even}
\end{array}\right.\\
& =
\left\{\begin{array}{ll}
\ol{\omega} (alt({\sigma |_{p^{\prime}_1}}), alt({\ul{k}})
alt({j^2_1} , \cdots,
  {j^2_{s-2r}}),g_{j^2_{s-2r}}) & \text{if $s$ is odd}\\
\omega (alt({\sigma |_{p^{\prime}_1}}), alt({\ul{k}})
alt({j^2_1} , \cdots,
  {j^2_{s-2r-1}}),g_{j^2_{s-2r}}) & \text{if $s$ is even}
\end{array}\right.\\
& \text{(applying (ii) and (iii))}\\
& = \; \lambda_{\sigma |_{p^{\prime}_1}}
((\ul{k},\ul{j}^2,\widetilde{\ul{k}}),s).\\
\end{align*}

\noindent
Third term: Note that
\begin{align*}
alt({i^1_{2r+1}} , \cdots , {i^1_{2r+n_2}}) \; = \; & \;
alt({j^2_1} , \cdots, {j^2_{n_2}})
\text{ (using (ii))}  \tag{v}\\
= \; & \; alt({j^2_{2n_2}} , \cdots,
{j^2_{n_2 +1}}) \text{ (since $alt({\ul{j}^2}) = e$)}\\
= \; & \; alt({j^1_{2r+1}} , \cdots , {j^1_{2r+n_2}})
\text{ (using (iv)).}
\end{align*}

Thus, for $2r+n_2+1 \leq s \leq 2n_1$,
\begin{align*}
& \; \lambda_{\sigma |_{p^{\prime}_1}} (\ul{i}^1, s)\\
= &
\left\{\begin{array}{ll}
\ol{\omega}(alt({\sigma |_{p^{\prime}_1}}) , alt({i^1_1} , \cdots
, {i^1_{2r}}) alt({i^1_{2r+1}} , \cdots, {i^1_{2r+n_2}}, \cdots
{i^1_s} , g_{i^1_s}) & \text{if $s$ is odd}\\
\omega(alt({\sigma |_{p^{\prime}_1}}) , alt({i^1_1} , \cdots
, {i^1_{2r}}) alt({i^1_{2r+1}} , \cdots, {i^1_{2r+n_2}}, \cdots
{i^1_{s-1}} , g_{i^1_s}) & \text{if $s$ is even}
\end{array}\right.\\
= &
\left\{\begin{array}{ll}
\ol{\omega}(alt({\sigma |_{p^{\prime}_1}}) , alt({j^1_1} , \cdots
, {j^1_s}) , g_{j^1_s}) & \text{if $s$ is odd}\\
\omega(alt({\sigma |_{p^{\prime}_1}}) , alt({j^1_1} , \cdots
, {j^1_{s-1}}) , g_{j^1_s}) & \text{if $s$ is even}
\end{array}\right.
\text{(using (i) and (v))}\\
 = & \;\lambda_{\sigma |_{p^{\prime}_1}} (\ul{j}^1, s).
\end{align*}

Combining the three terms, we get
\begin{align*}
\lambda_{\sigma |_{p^{\prime}_1}} (\ul{i}^1) & =  \left(
  \prod^{2r}_{s=1} \lambda_{\sigma |_{p^{\prime}_1}}
  (\ul{j}^1, s) \right)
\left( \prod^{2r+n_2}_{s=2r+1} \lambda_{\sigma |_{p^{\prime}_1}}
  ((\ul{k},\ul{j}^2,\widetilde{\ul{k}}), s) \right)
\left( \prod^{2n_1}_{s= 2r + n_2 + 1} \lambda_{\sigma |_{p^{\prime}_1}}
  (\ul{j}^1, s) \right)\\
& = \lambda_{\sigma |_{p^{\prime}_1}} (\ul{j}^1)
\left(
  \prod^{2r+n_2}_{s=2r+1} \ol{\lambda}_{\sigma |_{p^{\prime}_1}}
  (\ul{j}^1, s) \right)
\left( \prod^{2r+n_2}_{s=2r+1} \lambda_{\sigma |_{p^{\prime}_1}}
  ((\ul{k},\ul{j}^2,\widetilde{\ul{k}}), s) \right).
\end{align*}

Now, for $2r+1 \leq s \leq 2r+n_2$, we compute

\begin{align*}
& \; \ol{\lambda}_{\sigma |_{p^{\prime}_1}} (\ul{j}^1, s)\\
= & \left\{\begin{array}{ll}
\omega(alt({\sigma |_{p^{\prime}_1}}) , alt({j^1_1} , \cdots ,
j^1_{2r}, \cdots, {j^1_s}) , g_{j^1_s}) & \text{if $s$ is odd}\\
\ol{\omega}(alt({\sigma |_{p^{\prime}_1}}) , alt({j^1_1} , \cdots ,
j^1_{2r}, \cdots, {j^1_{s-1}}) , g_{j^1_s}) & \text{if $s$ is even}
\end{array}\right.\\
= & \left\{\begin{array}{ll}
\omega(alt({\sigma |_{p^{\prime}_1}}) , alt({\ul{k}})alt({j^2_{2n_2}},
 \cdots, {j^2_{2n_2+2r-s+1}}) , g_{j^2_{2n_2+2r-s+1}}) & \text{if
  $s$ is odd}\\
\ol{\omega}(alt({\sigma |_{p^{\prime}_1}}) , alt({\ul{k}})
alt({j^2_{2n_2}}, \cdots, {j^2_{2n_2+2r-s+2}}) , g_{j^2_{2n_2+2r-s+1}})
& \text{if $s$ is even}
\end{array}\right.\\
& \text{(using (iii) and (iv))}\\
= & \left\{\begin{array}{ll}
\omega(alt({\sigma |_{p^{\prime}_1}}) , alt({\ul{k}}) alt({j^2_1},
\cdots, {j^2_{2n_2+2r-s}}) , g_{j^2_{2n_2+2r-s+1}}) & \text{if
  $s$ is odd}\\
\ol{\omega}(alt({\sigma |_{p^{\prime}_1}}) , alt({\ul{k}})
alt({j^2_1}, \cdots, {j^2_{2n_2+2r-s+1}}) , g_{j^2_{2n_2+2r-s+1}}) & \text{if
$s$ is even}
\end{array}\right.\\
& \text{(since $alt({\ul{j}^2}) = e$)}\\
= & \lambda_{\sigma |_{p^{\prime}_1}}
  ((\ul{k},\ul{j}^2,\widetilde{\ul{k}}), 2n_2+4r-s+1).
\end{align*}

Hence, we obtain

\begin{align*}
\lambda_{\sigma |_{p^{\prime}_1}} (\ul{i}^1)
& = \lambda_{\sigma |_{p^{\prime}_1}} (\ul{j}^1)
\left( \prod^{2r+2n_2}_{s=2r+n_2+1} \lambda_{\sigma |_{p^{\prime}_1}}
  ((\ul{k},\ul{j}^2,\widetilde{\ul{k}}), s) \right)
\left( \prod^{2r+n_2}_{s=2r+1} \lambda_{\sigma |_{p^{\prime}_1}}
  ((\ul{k},\ul{j}^2,\widetilde{\ul{k}}), s) \right)\\
& = \lambda_{\sigma |_{p^{\prime}_1}} (\ul{j}^1) \;
\lambda_{\sigma |_{p^{\prime}_1}}
(\ul{k},\ul{j}^2,\widetilde{\ul{k}}) \; \text{ (by Remark
  \ref{leftinc2}).}
\end{align*}

We can now proceed with the proof of equation (\ref{discinc}):
\begin{align*}
\lambda_{\sigma |_{p^{\prime}_1}} (\ul{i}^1) \;
\lambda_{\tau |_{p_2}} (\ul{j}^2)  & = \lambda_{\sigma
  |_{p^{\prime}_1}} (\ul{j}^1) \;
\lambda_{\sigma |_{p^{\prime}_1}}
(\ul{k},\ul{j}^2,\widetilde{\ul{k}}) \; \lambda_{\ul{k}}
(\ul{j}^2)\\
& = \lambda_{\sigma
  |_{p^{\prime}_1}} (\ul{j}^1) \;
\; \lambda_{(\sigma
  |_{p^{\prime}_1}, \ul{k})} (\ul{j}^2),
\end{align*}

where we applied Lemma \ref{leftinc1} for the last identity.
This completes the proof that the action of tangles defined above
is compatible with composition of tangles (called
{\it naturality of composition} in \cite{J2}).
Hence, $P^{\langle g_i : i \in I \rangle, \omega}$ is a
planar algebra.

We will define next a {$*$-{\it structure} on $P^{\langle g_i :
i \in I \rangle, \omega}$}. Note that if
$\ul{i} \in I^{2n}$, then $alt(\ul{i}) = e$ if and only if
$alt(\widetilde{\ul{i}}) = e$. Extend the operation $\sim$ conjugate
linearly to define a $*$-structure on $P^{\langle g_i : i \in I \rangle,
\omega}_n$. Clearly, this is an involution. We need to check whether the
action of tangles preserves $*$, that is, $Z_{T^*} \circ (* \times
\cdots \times *) = * \circ Z_T$. It is enough to check this identity
when $T$ has no internal discs or closed loops, and
when $T$ is an elementary annular tangle.

If $T$ has no
internal discs or closed loops, then it is a Temperley-Lieb diagram and
hence $Z_T$ is the sum of all sequences of indices from $I$ which have
non-crossing matched pairings where the position of the pairings are
given by the numberings of the marked points on the boundary of $T$
which are connected by a string. Now, in the tangle $T^*$, the $m$-th
and the $n$-th
marked points are connected by a string if and only if the $m$-th and
the $n$-th marked points starting from the last point in $T$ reading
anticlockwise, are connected. So, $Z_{T^*}$ is indeed the sum of all
sequences featuring in the linear expansion of $Z_T$ in the reverse
order (that is, applying $\sim$).

If $T$ is an elementary annular tangle of capping
(resp. cap-inclusion) type with $m$-th and $(m+1)$-th marked points of
the internal (resp. external) disc being connected by a string, then
$T^*$ is also the same kind of elementary annular tangle but the `capping'
occurs at the $m$-th and $(m+1)$-th marked points of the internal
(resp. external) disc starting from the last point reading
anticlockwise. The identity easily follows from this observation.

If $T$ is an elementary tangle of left-inclusion type, then $T =
T^*$. The identity will then follow from the next lemma.

\begin{lemma}\label{lambda*}
If $\ul{j} \in I^m$ and $\ul{i} = (i_1, \cdots, i_{2n}) \in I^{2n}$ such
that $alt(\ul{i}) = e$, then $\lambda_{\ul{j}} (\widetilde{\ul{i}}) =
\ol{\lambda}_{\ul{j}} (\ul{i})$.
\end{lemma}

\noindent{\bf Proof:} The proof is an immediate consequence of
$alt(i_1, \cdots, i_s) = alt(i_{2n}, \cdots, i_{s+1})$,
which holds since $alt(\ul{i}) = e$. \qed

\medskip
Lastly, if $T$ is an elementary tangle of disc-inclusion type
given by Figure \ref{path}, then we need to show $\langle Z_{T^*}
(\widetilde{\ul{i}}, \widetilde{\ul{j}}), \widetilde{\ul{k}} \rangle =
\ol{\langle Z_T(\ul{i}, \ul{j}), \ul{k} \rangle}$ for $\ul{i}, \ul{k}
\in I^{2n_1}, \ul{j} \in I^{2n_2}$ with $alt(\ul{i}) =
alt(\ul{j}) = alt(\ul{k}) = e$. It is easy to verify that $\ul{i}$,
$\ul{j}$ and $\ul{k}$ define a state on $T$ if and only if
$\widetilde{\ul{i}}$, $\widetilde{\ul{j}}$ and $\widetilde{\ul{k}}$
define the same on $T^*$. If they fail to define a state, then both
sides are zero. If they define a state, then the scalars appearing on
both sides can be made equal by applying Lemma \ref{lambda*}.

\begin{remark}
If $\omega$ is trivial, that is, a coboundary, then $P^{\langle
g_i : i \in I \rangle, \omega}$ is isomorphic to example 2.7
in \cite{J2}.  Jones constructed this example of a planar algebra
by considering a certain subspace of the tensor planar algebra (TPA)
over the vector space with the indexing set $I$ as a basis.
He then showed that this
subspace is closed under the TPA-action of tangles, thus showing
that the subspace is indeed a planar algebra. One can view our planar
algebra $P^{\langle g_i : i \in I \rangle, \omega}$ as a subspace of
the TPA in an obvious way but the action of planar tangles induced
by TPA will not be the same as our action which involves the extra data of
a $3$-cocycle. It does not seem clear if
our planar algebra $P^{\langle g_i : i \in I \rangle, \omega}$
can be viewed as a planar subalgebra of the tensor planar algebra
over the vector space with basis $I$ if $\omega$ is nontrivial.
\end{remark}
\section{The planar algebra of the diagonal subfactor with
cocycle}\label{subfact}
In this section we will compute the relative commutants of the
diagonal subfactor associated to a $G$-kernel.
We will determine the filtered $*$-algebra structure,
Jones projections and the conditional
expectations. Note that some of this can already be found
in \cite{Bi1}, \cite{P2}, \cite{J2}. The main point here is that we are
able to choose an appropriate basis of the higher relative commutants
such that the action of planar tangles on these allows us to identify
this concrete planar algebra with the abstract one defined in the
previous section.
\noindent \vskip 1em
Let $N$ be a II$_1$ factor, $I$ be a finite set and for $i \in I$,
choose $\theta_i \in \Aut N$. Set $M = M_I \otimes N$ where $M_I$
denotes the matrix algebra whose rows and columns are indexed by $I$.
As in the introduction, consider the {\it diagonal subfactor}
$N \hookrightarrow M$ given by
\[N \ni x \longmapsto \sum_{i \in I} E_{i,i} \otimes \theta_i (x) \in
M\]
that is, an element $x$ of
$N$ sits in $M$ as a diagonal matrix whose $i$th diagonal is $\theta_i
(x)$. If we
have another collection of automorphisms $\theta^{\prime}_i \in \Aut N$
for $i \in I$ such that $\theta_i =
\theta^{\prime}_i \mod \Int N$, for all $i \in I$ (up to a permutation
of $I$), then the associated diagonal subfactors
are isomorphic. It is therefore natural to associate a diagonal
subfactor to a collection of elements
$g_i \in \Out N$, $i \in I$. We consider the
subgroup $G = \langle g_i : i \in I \rangle$ of $\Out N$, which can
be viewed as a $G$-kernel in the obvious way, and choose a lift
\[ G \ni g \stackrel{\alpha}{\longmapsto} \alpha_g \in \Aut N\]
such that $\alpha_e = id_N$. ~
Set $\alpha_i = \alpha_{g_i}$ for all $i \in I$. Consider the diagonal
subfactor $N \subset M = M_I \otimes N$ where
the $i$th diagonal entry of an element of $N$ viewed in $M$ is twisted
by the action of $\alpha_{i}$ . The index of this subfactor is $|I|^2$.
Set $M_n = M_{I^{n+1}} \otimes N$ for $n \geq 0$. We will often
identify $E_{\ul{i},\ul{j}} \otimes E_{\ul{k},\ul{l}} \in M_{I^m}
\otimes M_{I^n}$ with
$E_{(\ul{i},\ul{k}),(\ul{j},\ul{l})} \in M_{I^{m+n}}$ for $\ul{i},
\ul{j} \in I^m$ and $\ul{k}, \ul{l} \in I^n$. Now, $M_{n-1}$ is
included in $M_n$ in the following way:
\[M_{n-1} = M_{I^n} \otimes N \ni E_{\ul{i},\ul{j}} \otimes x
\longmapsto E_{\ul{i},\ul{j}} \otimes \psi_n (x) \in M_{I^n} \otimes
M_I \otimes N = M_{I^{n+1}} \otimes N = M_n\]
for all $\ul{i}, \ul{j} \in I^n$ and $x \in N$ where $\psi_n : N
\rightarrow M_I \otimes N$ is defined as:
\[
\psi_n (x) =
\left\{\begin{array}{ll}
\sum_{k \in I} E_{k,k} \otimes \alpha_{k} (x) & \text{ if $n$ is
  even}\\
\sum_{k \in I} E_{k,k} \otimes \alpha^{-1}_{k} (x) & \text{ if $n$
  is odd}
\end{array}\right.
\]
It is easy to check (see \cite{P2}, \cite{Bi1}) that
$N \subset M \subset M_1 \subset M_2 \subset \cdots$ is isomorphic
to the Jones tower of II$_1$ factors associated to
$N \subset M$ where the Jones projections and
conditional expectations are given by:

\begin{align*}
e_n & =  \; |I|^{-1} \sum_{\ul{k} \in I^{n-1}, i,j \in I}
E_{(\ul{k},i,i),(\ul{k},j,j)} \otimes 1 \in M_n\\
\E^{M_n}_{M_{n-1}} (E_{(\ul{i},k),(\ul{j},l)} \otimes x) & =
\delta_{k,l} \; |I|^{-1} \; E_{\ul{i},\ul{j}} \otimes
\alpha^{(-1)^{n-1}}_{k} (x)
\end{align*}
for all $\ul{i}, \ul{j} \in I^n$, $k,l \in I$, $x \in
N$. Moreover, $\left\{\sqrt{|I|} (E_{i,j} \otimes 1) : i,j
\in I \right\}$ forms a Pimsner-Popa basis of $M$ over
$N$. This basis will be used to write the conditional expectation of
commutant of $N$ onto the commutant of $M$ (see \cite{Bi2}).

To find the relative commutant $N^\prime \cap M_{n-1}$, first note
that $N$ is included in $M_{n-1}$ by the following map:
\[N \ni x \longmapsto \sum_{\ul{i} \in I^n} E_{\ul{i},\ul{i}}
\otimes alt^{-1}_\alpha (\ul{i}) (x) \in M_{n-1}\]
where $alt_\alpha (\ul{i}) = \alpha^{-1}_{i_1}\alpha_{i_2} \cdots
\alpha^{(-1)^n}_{i_n} \in \Aut N$ for $\ul{i} = (i_1, \cdots, i_n) \in
I^n$. Now, if $\sum_{\ul{i},\ul{j} \in I^n} x_{\ul{i},\ul{j}}
(E_{\ul{i},\ul{j}} \otimes 1) \in N^\prime \cap M_{n-1}$, then
\begin{align*}
\sum_{\ul{i},\ul{j} \in I^n} x_{\ul{i},\ul{j}}
(E_{\ul{i},\ul{j}} \otimes 1) \; y  & = y \sum_{\ul{i},\ul{j} \in I^n}
x_{\ul{i},\ul{j}} (E_{\ul{i},\ul{j}} \otimes 1)  \text{ for all $y
  \in N$}\\
\Leftrightarrow  x_{\ul{i},\ul{j}} \left(alt_\alpha (\ul{i})
alt^{-1}_\alpha (\ul{j}) \right) (y) & = y \; x_{\ul{i},\ul{j}}
\text{ for all $y \in N$, $\ul{i}, \ul{j} \in I^n$}\\
\Leftrightarrow  x_{\ul{i},\ul{j}} alt_\alpha (\ul{i} ,
\widetilde{\ul{j}}) (y) & = y \; x_{\ul{i},\ul{j}}
\text{ for all $y \in N$, $\ul{i}, \ul{j} \in I^n$}
\end{align*}

Thus, for $\ul{i}, \ul{j} \in I^n$, if $x_{\ul{i},\ul{j}} \neq 0$, then
$x_0 =\frac{x_{\ul{i},\ul{j}}}{\|x_{\ul{i},\ul{j}}\|} \in \mathcal{U}
(N)$ and $\Ad {x_0} \circ alt_\alpha (\ul{i}, \widetilde{\ul{j}}) =
id_N$ which implies $alt(\ul{i},\widetilde{\ul{j}}) = e$. Similarly, if there
exist $\ul{i}, \ul{j} \in I^n$ such that $alt(\ul{i},\widetilde{\ul{j}}) = e$,
then $u(E_{\ul{i},\ul{j}} \otimes 1) \in N^\prime \cap M_{n-1}$ where
$u \in \mathcal{U} (N)$ satisfies $\Ad u \circ alt_\alpha (\ul{i},
\widetilde{\ul{j}}) = id_N$. Thus,
\[N^\prime \cap M_{n-1} = span\left\{u(E_{\ul{i},\ul{j}} \otimes 1)
\in M_{n-1}
\left|\begin{array}{c}
\ul{i}, \ul{j} \in I^n \text{ and }  u \in \mathcal{U} (N)\\
\text{s.t. }\Ad u \circ alt_\alpha(\ul{i},\widetilde{\ul{j}}) = id_N
\end{array}\right.
\right\}.\]

The elements in this set do not yet form a basis since
the unitary $u$ can be modified by a scalar of absolute value $1$.
To get a good basis of $N^\prime \cap M_{n-1}$ we need to choose $u$
in such a way that the planar algebra associated
to $N \subset M$ can easily be identified with the abstract one defined
in Section \ref{plnalg}.

We now digress a little bit to set up some notations. Let $u: G \times
G \rightarrow \mathcal{U} (N)$ be a unitary defined by

\[\alpha_{g_1} \alpha_{g_2} = \Ad {u(g_1,g_2)} \circ \alpha_{g_1 g_2}
\; \text{ for all } g_1, g_2 \in G\]
such that $u(g_1,g_2) = 1$ whenever either of $g_1$ or $g_2$ is
$e$. For $\ul{i} = (i_1, \cdots, i_n) \in I^n$, define
\[v_m(\ul{i}) =
\left\{ \begin{array}{ll}
u^*(alt(i_1, \cdots, i_{m}),g_{i_m}) & \text{if $m$ is odd}\\
u(alt(i_1, \cdots, i_{m-1}),g_{i_m}) & \text{if $m$ is even}
\end{array}
\right.\]
and set $v(\ul{i}) = v_1(\ul{i}) \cdots v_n(\ul{i})$. Next, we prove a
several useful lemmas involving $v$.
The first lemma motivates our choice of the basis.
\begin{lemma}\label{unichoice}
$alt_\alpha (\ul{i}) =\Ad{v(\ul{i})} \circ \alpha_{alt(\ul{i})}$ for
  all $\ul{i} \in I^n$.
\end{lemma}
\noindent{\bf Proof:} Using the definition of $u$, note that for all
$m \geq 1$,
\[\Ad {v_m(\ul{i})} =
\left\{\begin{array}{ll}
\alpha_{alt(i_1, \cdots, i_{m-1})} \alpha^{-1}_{i_m}
\alpha^{-1}_{alt(i_1, \cdots, i_m)} & \text{if $m$ is odd}\\
\alpha_{alt(i_1, \cdots, i_{m-1})} \alpha_{i_m} \alpha^{-1}_{alt(i_1,
  \cdots, i_m)} & \text{if $m$ is even}
\end{array}\right.\]
Hence
\[\Ad{v(\ul{i})} =\Ad{v_1(\ul{i})} \cdots \Ad{v_n(\ul{i})} = \alpha_e
\alpha^{-1}_{i_1} \alpha_{i_2} \cdots \alpha^{(-1)^n}_{i_n}
\alpha^{-1}_{alt(i_1, \cdots, i_n)} = alt_\alpha (\ul{i})
\alpha^{-1}_{alt(\ul{i})}.\]\qed

\begin{lemma}\label{v*}
For any $\ul{k} = (k_1, \cdots, k_{2n}) \in I^{2n}$ such that
$alt(\ul{k}) = e$, we have $v(\widetilde{\ul{k}}) v(\ul{k}) = 1$.
\end{lemma}

\noindent{\bf Proof:} First we expand $v(\widetilde{\ul{k}})$ and
$v(\ul{k})$ into products of unitaries arising from the definition of
$v$, and then we consider that the product of the $p$-th unitary of
$v(\widetilde{\ul{k}})$ from the right and $p$-th unitary of
$v(\ul{k})$ from the left, that is,
\begin{align*}
v_{2n-p+1}(\widetilde{\ul{k}}) v_p(\ul{k}) & =
\left\{\begin{array}{ll}
u(alt(k_{2n}, \cdots, k_{p+1}), k_p) \; u^*(alt(k_1, \cdots, k_p),
k_p) & \text{if $p$ is odd}\\
u^*(alt(k_{2n}, \cdots, k_p), k_p) \; u(alt(k_1, \cdots, k_{p-1}),
k_p) & \text{if $p$ is even}
\end{array}\right.\\
& =
\left\{\begin{array}{ll}
u(alt(k_1, \cdots, k_p), k_p) \; u^*(alt(k_1, \cdots, k_p),
k_p) & \text{if $p$ is odd}\\
u^*(alt(k_1, \cdots, k_{p-1}), k_p) \; u(alt(k_1, \cdots, k_{p-1}),
k_p) & \text{if $p$ is even}
\end{array}\right.\\
& = 1,
\end{align*}
for $1 \leq p \leq 2n$.\qed

\medskip
\begin{lemma}\label{vmult}
For $\ul{i} = (i_1, \cdots, i_n), \; \ul{j} = (j_1, \cdots, j_n), \;
\ul{k} = (k_1, \cdots, k_n) \in I^n$ such that $alt(\ul{i},
\widetilde{\ul{j}}) = e = alt(\ul{j}, \widetilde{\ul{k}})$, we have
$v(\ul{i}, \widetilde{\ul{j}}) \; v(\ul{j}, \widetilde{\ul{k}}) =
v(\ul{i}, \widetilde{\ul{k}})$.
\end{lemma}

\noindent{\bf Proof:} Using an argument similar to the proof of Lemma
\ref{v*}, one can prove that the product of the $p$-th unitary of
$v(\ul{i} , \widetilde{\ul{j}})$ from the right and the $p$-th unitary
of $v(\ul{j}, \widetilde{\ul{k}})$ from the left, is $1$ for $1 \leq p
\leq n$. Again, for $n+1 \leq p \leq 2n$,
\begin{align*}
v_p(\ul{j}, \widetilde{\ul{k}}) & =
\left\{\begin{array}{ll}
u^*(alt(\ul{j}) k^{(-1)^{n+1}}_{2n} \cdots k^{-1}_{2n-p+1},
k_{2n-p+1}) & \text{if $p$ is odd}\\
u(alt(\ul{j})k^{(-1)^{n+1}}_{2n} \cdots k^{-1}_{2n-p+2}, k_{2n-p+1}) &
\text{if $p$ is even}
\end{array}\right.\\
& =
\left\{\begin{array}{ll}
u^*(alt(\ul{i}) k^{(-1)^{n+1}}_{2n} \cdots k^{-1}_{2n-p+1},
k_{2n-p+1}) & \text{if $p$ is odd}\\
u(alt(\ul{i})k^{(-1)^{n+1}}_{2n} \cdots k^{-1}_{2n-p+2}, k_{2n-p+1}) &
\text{if $p$ is even}
\end{array}\right.
\text{($alt(\ul{i}, \widetilde{\ul{j}}) = e$)}\\
& = v_p(\ul{i}, \widetilde{\ul{k}}).
\end{align*}
Thus,
\begin{align*}
v(\ul{i}, \widetilde{\ul{j}}) \; v(\ul{j}, \widetilde{\ul{k}}) & =
v_1(\ul{i}, \widetilde{\ul{j}}) \cdots v_n(\ul{i}, \widetilde{\ul{j}})
v_{n+1}(\ul{j}, \widetilde{\ul{k}}) \cdots v_{2n}(\ul{j},
\widetilde{\ul{k}})\\
& = v_1(\ul{i}, \widetilde{\ul{k}}) \cdots v_n(\ul{i},
\widetilde{\ul{k}}) v_{n+1}(\ul{i}, \widetilde{\ul{k}}) \cdots
v_{2n}(\ul{i}, \widetilde{\ul{k}}) = v(\ul{i}, \widetilde{\ul{k}})
\end{align*}
\qed\\

By applying Lemma \ref{unichoice}, we see that the set
$\left\{v^*(\ul{i},\widetilde{\ul{j}}) (E_{\ul{i},\ul{j}} \otimes 1)
\left|\begin{array}{l}
\ul{i},\ul{j} \in I^n \text{ s.t.}\\
alt(\ul{i},\widetilde{\ul{j}}) = e
\end{array}\right.
\right\}$
is a basis for $N^\prime \cap M_{n-1}$. We will use this basis to
establish an isomorphism between the planar algebra associated to $N
\subset M$ and the abstract planar algebra defined in
Section \ref{plnalg}. Let $\omega : G \times G \times G
\rightarrow \mathbb T$ be the $3$-cocycle associated to $G$, that is,
for all $g_1, g_2, g_3 \in G$ we have
\[u(g_1,g_2) u(g_1g_2,g_3) = \omega(g_1,g_2,g_3) \alpha_{g_1}
(u(g_2,g_3)) u(g_1,g_2g_3).\]
As before, this is a consequence of associativity of the multiplication in
$G$.

We prove next a useful lemma relating $v$ and $\omega$.

\begin{lemma}\label{vleftcondexp}
If $i \in I$ and $\ul{k}=(k_1, \cdots, k_{2n}) \in I^{2n}$ such that
$alt(\ul{k}) = e$ and $k_1 = k_{2n}$, then $alt_\alpha (i,k_1)
(v(\ul{k})) = \ol{\lambda}_{(i,k_1)} (\ul{k}) \; v(i, k_2, \cdots,
k_{2n-1}, i)$.
\end{lemma}

\noindent{\bf Proof:} We expand $alt_\alpha (i,k_1) (v(\ul{k}))$ as a
product of unitaries and work with them separately. For $1 \leq p \leq
n$, we have\\
$\begin{array}{lll}
\text{(i)} & & alt_\alpha (i,k_1) (v_{2p-1}(\ul{k}))\\
& = & \left(\alpha^{-1}_{g_i} \alpha_{k_1}\right) (u^*(alt(k_1,
\cdots, k_{2p-1}), g_{k_{2p-1}}))\\
& = & \left(\Ad{u^*(g^{-1}_i,g_i) u(g^{-1}_i, g_{k_1})} \circ
\alpha_{g^{-1}_i g_{k_1}} \right) (u^*(alt(k_1, \cdots, k_{2p-1}),
g_{k_{2p-1}}))\\
& & \text{(using definition of $u$)}\\
& = & \; \omega(alt(i,k_1), alt(k_1, \cdots, k_{2p-1}), g_{k_{2p-1}})\\
& & \Ad{u^*(g^{-1}_i,g_i) u(g^{-1}_i, g_{k_1})}
\left(
\begin{array}{l}
u(alt(i,k_1), alt(k_1, \cdots, k_{2p-2}, k_{2p-1}, k_{2p-1}))\\
\cdot \; u^*(alt(i, k_1, k_1, k_2, \cdots, k_{2p-1}), g_{k_{2p-1}})\\
\cdot \; u^*(alt(i,k_1), alt(k_1, \cdots, k_{2p-1}))
\end{array}
\right)\\
& & \text{(using the definition of $\omega$)}
\end{array}$\\
$\begin{array}{lll}
\text{(ii)} & & alt_\alpha (i,k_1) (v_{2p}(\ul{k}))\\
& = & \left(\alpha^{-1}_{g_i} \alpha_{k_{i_1}}\right) (u(alt(k_1,
\cdots, k_{2p-1}), g_{k_{2p}}))\\
& = & \left(\Ad{u^*(g^{-1}_i,g_i) u(g^{-1}_i, g_{k_1})} \circ
\alpha_{g^{-1}_i g_{k_1}} \right) (u(alt(k_1, \cdots, k_{2p-1}),
g_{k_{2p}}))\\
& & \text{(using definition of $u$)}\\
& = & \; \ol{\omega}(alt(i,k_1), alt(k_1, \cdots, k_{2p-1}),
g_{k_{2p}})\\
& & \Ad{u^*(g^{-1}_i,g_i) u(g^{-1}_i, g_{k_1})}
\left(
\begin{array}{l}
u(alt(i,k_1), alt(k_1, \cdots, k_{2p-1}))\\
\cdot \; u(alt(i, k_1, k_1, k_2, \cdots, k_{2p-1}), g_{k_{2p}})\\
\cdot \; u^*(alt(i,k_1), alt(k_1, \cdots, k_{2p}))
\end{array}
\right)\\
& & \text{(using the definition of $\omega$)}
\end{array}$\\
Multiplying (i) and (ii), we get,
\begin{align*}
& alt_\alpha (i,k_1) ((v_{2p-1}(\ul{k}))\;(v_{2p}(\ul{k})))\\
= & \ol{\lambda}_{(i,k_1)} (\ul{k}, 2p-1) \; \ol{\lambda}_{(i,k_1)}
(\ul{k}, 2p)\\
& \Ad{u^*(g^{-1}_i,g_i) u(g^{-1}_i, g_{k_1})}
\left(
\begin{array}{l}
u(alt(i,k_1), alt(k_1, \cdots, k_{2p-2}, k_{2p-1}, k_{2p-1}))\\
\cdot \; u^*(alt(i, k_1, k_1, k_2, \cdots, k_{2p-1}), g_{k_{2p-1}})\\
\cdot \; u(alt(i, k_1, k_1, k_2, \cdots, k_{2p-1}), g_{k_{2p}})\\
\cdot \; u^*(alt(i,k_1), alt(k_1, \cdots, k_{2p}))
\end{array}
\right)\\
= & \ol{\lambda}_{(i,k_1)} (\ul{k}, 2p-1) \; \ol{\lambda}_{(i,k_1)}
(\ul{k}, 2p)\\
&
\left\{\begin{array}{ll}
\Ad{u^*(g^{-1}_i,g_i) u(g^{-1}_i, g_{k_1})}
\left(
\begin{array}{l}
u(alt(i,k_1), e) \; u^*(g^{-1}_i, g_{k_1})\\
\cdot \; v_2(i, k_2, \cdots, k_{2n-1}, i)\\
\cdot \; u^*(alt(i,k_1), alt(k_1, k_2))
\end{array}
\right) & \text{if $p = 1$}\\
&\\
\Ad{u^*(g^{-1}_i,g_i) u(g^{-1}_i, g_{k_1})}
\left(
\begin{array}{l}
u(alt(i,k_1), alt(k_1, \cdots, k_{2p-2}))\\
\cdot \; v_{2p-1}(i, k_2, \cdots, k_{2n-1}, i)\\
\cdot \; v_{2p}(i, k_2, \cdots, k_{2n-1}, i)\\
\cdot \; u^*(alt(i,k_1), alt(k_1, \cdots, k_{2p}))
\end{array}
\right) & \text{if $2 \leq p \leq n-1$}\\
&\\
\Ad{u^*(g^{-1}_i,g_i) u(g^{-1}_i, g_{k_1})}
\left(
\begin{array}{l}
u(alt(i,k_1), alt(k_1, \cdots, k_{2n-2}))\\
\cdot \; v_{2n-1}(i, k_2, \cdots, k_{2n-1}, i)\\
\cdot \; u(g^{-1}_i, g_{k_1}) \; u^*(alt(i,k_1), e)
\end{array}
\right) & \text{if $p=n$}
\end{array}\right.\\
\overset{\text{\rm def}} = & W_p.
\end{align*}
Thus,
\begin{align*}
& alt_\alpha (i,k_1) (v(\ul{k})) = W_1 \cdots W_n =
  \ol{\lambda}_{(i,k_1)} (\ul{k}) \cdot\\
& \Ad{u^*(g^{-1}_i,g_i) u(g^{-1}_i, g_{k_1})}
\left(
\begin{array}{l}
u^*(g^{-1}_i, g_{k_1})\\
\cdot \; v_2(i, k_2, \cdots, k_{2n-1}, i) \cdots v_{2n-1}(i, k_2,
\cdots, k_{2n-1}, i)\\
\cdot \; u(g^{-1}_i, g_{k_1})
\end{array}
\right)\\
= & \; \ol{\lambda}_{(i,k_1)} (\ul{k}) u^*(g^{-1}_i,g_i) v_2(i, k_2,
\cdots, k_{2n-1}, i) \cdots v_{2n-1}(i, k_2, \cdots, k_{2n-1}, i)
u(g^{-1}_i,g_i)\\
= & \; \ol{\lambda}_{(i,k_1)} (\ul{k}) \; v(i, k_2, \cdots, k_{2n-1}, i).
\end{align*}\qed

Let us recall the following well-known fact about isomporphisms of
two planar algebras which will be used in the next theorem (\cite{J2}).
Let $P^1$ and $P^2$ be two planar algebras. Then
$P^1 \cong P^2$ (as planar algebras) if and only if there exist
a vector space isomorphism
$\psi : P^1_n \rightarrow P^2_n$ such that:
\begin{itemize}
\item[(i)] $\psi$ preserves the filtered algebra structure,
\item[(ii)] $\psi$ preserves the actions of all Jones
  projection tangles and the (two types of) conditional expectation
  tangles.
\end{itemize}

If $P^1$ and $P^2$ are $*$-planar algebras, then we require $\psi$ to
be $*$-preserving.

\begin{theorem}\label{mainiso}
The planar algebra $P^{sf}$ associated to the diagonal subfactor
obtained from a
II$_1$ factor $N$ and a finite collection of automorphisms $\alpha_i
\in \Aut N$ for $i \in I$, is isomorphic to $P^{\langle g_i:i\in
  I\rangle, \omega}$ where $g_i = [\alpha_i] \in \Out N$ for
all $i \in I$, and $\omega$ is the normalized $3$-cocycle associated
to $G= \langle g_i:i\in I\rangle \subseteq \Out N$ as above. 
\end{theorem}
\noindent{\bf Proof:} Let $G = \langle g_i : i \in I\rangle$ and
without loss of generality, let us assume that $\alpha$ is a lift of
$G$ such that $\alpha_i = \alpha_{g_i}$. By \cite{J2} we have
$P^{sf}_n = N^\prime \cap M_{n-1}$ for all $n \geq
0$.

Define the map $\phi:P^{sf} \rightarrow P \overset{\text{\rm def}}=
P^{\langle g_i:i\in
I\rangle, \omega}$ by first defining it on basis elements as
$\phi(v^*(\ul{i},\widetilde{\ul{j}}) (E_{\ul{i},\ul{j}} \otimes 1)) =
(\ul{i},\widetilde{\ul{j}})$ for all $\ul{i}, \ul{j} \in I^n$ such
that $alt(\ul{i},\widetilde{\ul{j}}) = e$, and then extending it
linearly. Clearly, $\phi$ is a vector space isomorphism. We will
show that $\phi$ is $*$-planar algebra
isomorphism. We make first the following observation:

For $i \in I$, $\ul{i} = (i_1, \cdots, i_n) \in
I^n$ and $0 \leq s \leq n$, we have the identity $v(\ul{i}) = v(i_1,
\cdots, i_s, i, i, i_{s+1}, \cdots, i_n)$. The proof is similar to
that of Lemma \ref{pair1}. Thus, if
$\ul{i}$ is a sequence of indices with non-crossing matched pairings,
then using this identity several times to reduce all consecutive matched
pairings, we get $v(\ul{i}) = 1$.

We show now that $\phi$ is indeed a planar algebra isomorphism following
the remark just before the theorem.

\medskip
\noindent(a) {\it $\phi$ is unital}: Since $1_{P^{sf}_n} = \sum_{\ul{i}
  \in I^n} (E_{\ul{i}, \ul{i}} \otimes 1)$ and $v(\ul{i},
  \widetilde{\ul{i}}) = 1$ by the above observation, we get $\phi
  (1_{P^{sf}_n}) = \sum_{\ul{i} \in I^n} (\ul{i}, \widetilde{\ul{i}})
  = 1_{P_n}$.

\medskip
\noindent(b) {\it $\phi$ preserves Jones projection tangles}: By
Theorem 4.2.1 in \cite{J2}, the $n$-th Jones projection tangle
$E_n$ acts as $Z^{P^{sf}}_{E_n} = |I| e_n = \sum_{\ul{k} \in I^{n-1},
i,j \in I} (E_{(\ul{k},i,i),(\ul{k},j,j)} \otimes 1)$. Since $(\ul{k},
i, i , j, j, \widetilde{\ul{k}})$ is a sequence of indices with
non-crossing matched pairings, by the above note $v(\ul{k}, i, i , j,
j, \widetilde{\ul{k}}) = 1$. So, $\phi(Z^{P^{sf}}_{E_n}) =
\sum_{\ul{k} \in I^{n-1}, i,j \in I} (\ul{k}, i, i , j, j,
\widetilde{\ul{k}}) = Z^P_{E_n}$.

\medskip
\noindent(c) {\it $\phi$ preserves the action of conditional
  expectation tangle}: By Theorem 4.2.1 in \cite{J2}, the action
  of conditional expectation tangle $\mathcal{E}^{n+1}_n$ is given by
$Z^{P^{sf}}_{\mathcal{E}^{n+1}_n} = |I| \; \E^{M_n}_{M_{n-1}}
|_{N^\prime \cap M_n}$. We compute
\begin{align*}
Z^{P^{sf}}_{\mathcal{E}^{n+1}_n} (v^* (\ul{i}, k, l,
\widetilde{\ul{j}}) (E_{(\ul{i},k),(\ul{j},l)} \otimes 1)) & =
v^*(\ul{i}, k, l, \widetilde{\ul{j}}) \E^{M_n}_{M_{n-1}}
(E_{(\ul{i},k),(\ul{j},l)} \otimes 1)\\
& = \delta_{k,l} \; v^* (\ul{i}, k, k, \widetilde{\ul{j}})
(E_{\ul{i},\ul{j}} \otimes 1)\\
& = \delta_{k,l} \; v^* (\ul{i}, \widetilde{\ul{j}})
(E_{\ul{i},\ul{j}} \otimes 1) \stackrel{\phi}{\longmapsto} \delta_{k,l}
\; (\ul{i}, \widetilde{\ul{j}}) \in P_n
\end{align*}
for all $\ul{i}, \ul{j} \in I^n$, $k,l \in I$ such that $alt(\ul{i},
k, l, \widetilde{\ul{j}}) = e$. From the action of tangles defined in
Section \ref{plnalg}, it is easy to check $Z^P_{\mathcal{E}^{n+1}_n}
(\ul{i}, k, l, \widetilde{\ul{j}}) = \delta_{k,l} \; (\ul{i},
\widetilde{\ul{j}})$.

\medskip
\noindent(d) {\it $\phi$ preserves $*$}: Applying $*$ on a basis
element of $P^{sf}_n$, we get
\begin{align*}
\left(v^* (\ul{i}, \widetilde{\ul{j}}) (E_{\ul{i},\ul{j}} \otimes
1)\right)^* & = (E_{\ul{j},\ul{i}} \otimes 1) v(\ul{i},
\widetilde{\ul{j}}) = \left(alt_\alpha (\ul{j}) alt^{-1}_\alpha
(\ul{i}) \right) (v(\ul{i}, \widetilde{\ul{j}})) (E_{\ul{j},\ul{i}}
\otimes 1)\\
& = alt_\alpha (\ul{j}, \widetilde{\ul{i}}) (v(\ul{i},
\widetilde{\ul{j}})) (E_{\ul{j},\ul{i}} \otimes 1)\\
& = v(\ul{j}, \widetilde{\ul{i}}) \; v(\ul{i}, \widetilde{\ul{j}}) \;
v^*(\ul{j}, \widetilde{\ul{i}}) (E_{\ul{j},\ul{i}} \otimes 1)\\
& = v^*(\ul{j}, \widetilde{\ul{i}}) (E_{\ul{j},\ul{i}} \otimes 1)
\stackrel{\phi}{\longmapsto} (\ul{j}, \widetilde{\ul{i}}) = (\ul{i},
\widetilde{\ul{j}})^*
\end{align*}
for $\ul{i}, \ul{j} \in I^n$ such that $alt(\ul{i},\widetilde{\ul{j}})
= e$. Note that we used Lemma \ref{unichoice} for the third equality and
Lemma \ref{v*} for the fourth one.

\medskip
\noindent(e) {\it $\phi$ preserves multiplication}: Suppose $\ul{i},
\ul{j}, \ul{k}, \ul{l} \in I^n$ such that $alt(\ul{i},
\widetilde{\ul{j}}) = e = alt(\ul{k}, \widetilde{\ul{l}})$. Then,
\begin{align*}
& \left(v^*(\ul{i}, \widetilde{\ul{j}}) (E_{\ul{i},\ul{j}} \otimes 1)
\right) \cdot \left(v^*(\ul{k}, \widetilde{\ul{l}}) (E_{\ul{k},\ul{l}}
\otimes 1) \right)\\
= & \; v^*(\ul{i}, \widetilde{\ul{j}}) \left(alt_\alpha (\ul{i},
\widetilde{\ul{j}}) \right) (v^*(\ul{k}, \widetilde{\ul{l}}))
(E_{\ul{i},\ul{j}} \otimes 1) (E_{\ul{k},\ul{l}} \otimes 1)\\
= & \; v^*(\ul{i}, \widetilde{\ul{j}})  v(\ul{i}, \widetilde{\ul{j}})
v^*(\ul{k}, \widetilde{\ul{l}})) v^*(\ul{i}, \widetilde{\ul{j}})
\delta_{\ul{j} , \ul{k}} \; (E_{\ul{i},\ul{l}} \otimes 1) \; \text{
(using Lemma \ref{unichoice})}\\
= & \delta_{\ul{j} , \ul{k}} \left(v(\ul{i}, \widetilde{\ul{j}})
v(\ul{j}, \widetilde{\ul{l}})) \right)^* (E_{\ul{i},\ul{l}} \otimes 1)
= \delta_{\ul{j} , \ul{k}} \; v^*(\ul{i}, \widetilde{\ul{l}})
(E_{\ul{i},\ul{l}} \otimes 1) \; \text{ (using Lemma \ref{vmult})}
\end{align*}
On the other hand, one can easily deduce from the action of the
multiplication tangle in $P$ that $(\ul{i}, \widetilde{\ul{j}}) \cdot
(\ul{k}, \widetilde{\ul{l}}) = \delta_{\ul{j} , \ul{k}} \; (\ul{i},
\widetilde{\ul{l}})$.

\medskip
\noindent(f) {\it $\phi$ preserves the action of the left conditional
expectation tangle}: By Theorem 4.2.1 in \cite{J2}, the action of
the left conditional expectation tangle $\mathcal{E}^\prime_n$ is given by
$Z^{P^{sf}}_{\mathcal{E}^\prime_n} = |I| \; \E^{N^\prime \cap
M_{n-1}}_{M^\prime \cap M_{n-1}}$. Using the basis of $M$ over
$N$ mentioned before, the conditional expectation onto $M^\prime \cap
N_{n-1}$ can be expressed as (see \cite{Bi2})
\begin{align*}
\E^{N^\prime \cap M_{n-1}}_{M^\prime \cap M_{n-1}} (x) & =
|I|^{-2} \sum_{i,j \in I} \left(\sqrt{|I|} (E_{i,j} \otimes 1)\right) x
\left(\sqrt{|I|} (E_{j,i} \otimes 1)\right)\\
& = |I|^{-1} \sum_{i,j \in I} (E_{i,j} \otimes 1) x (E_{j,i}
\otimes 1)
\end{align*}
for $x \in N^\prime \cap M_{n-1}$. Hence, for $\ul{i} = (i_1, \cdots,
i_{n-1})$, $\ul{j} = (j_1, \cdots, j_{n-1}) \in I^{n-1}$ and $k, l \in
I$ such that $alt(k, \ul{i}, \widetilde{\ul{j}}, l) = e$, we have
\begin{align*}
& \; Z^{P^{sf}}_{\mathcal{E}^\prime_n} \left(v^*(k, \ul{i},
\widetilde{\ul{j}}, l) (E_{(k,\ul{i}),(l,\ul{j})} \otimes 1)\right)\\
= & \sum_{i,j \in I} (E_{i,j} \otimes 1) v^*(k, \ul{i},
\widetilde{\ul{j}}, l) (E_{(k,\ul{i}),(l,\ul{j})} \otimes 1) (E_{j,i}
\otimes 1)\\
= & \sum_{i,j \in I}  alt_\alpha (i,j) (v^*(k, \ul{i},
\widetilde{\ul{j}}, l)) (E_{i,j} \otimes 1) (E_{(k,\ul{i}),(l,\ul{j})}
\otimes 1) (E_{j,i} \otimes 1)\\
= & \; \delta_{k,l} \sum_{i \in I}  alt_\alpha (i,k) (v^*(k, \ul{i},
\widetilde{\ul{j}}, k)) (E_{(i,\ul{i}),(i,\ul{j})}
\otimes 1)\\
= & \; \delta_{k,l} \sum_{i \in I} \lambda_{i,k} (k, \ul{i}, \ul{j},
k) v^*(i, \ul{i}, \widetilde{\ul{j}}, i) (E_{(i,\ul{i}),(i,\ul{j})}
\otimes 1) \; \text{ (using Lemma \ref{vleftcondexp})}
\end{align*}
On the other hand, one can easily check that the action of
$\mathcal{E}^\prime_n$ is given by
\[Z^P_{\mathcal{E}^\prime_n} (k, \ul{i}, \widetilde{\ul{j}}, l) =
  \delta_{k,l} \sum_{i \in I} \lambda_{i,k} (k, \ul{i}, \ul{j}, k) \;
  (i, \ul{i}, \widetilde{\ul{j}}, i).\]\qed

\begin{corollary}
Given a group $G$ generated by a finite collection $g_i$ for $i \in I$
and given a normalized $3$-cocycle $\omega \in Z^3 (G,\mathbb T)$,
there exists a hyperfinite subfactor whose associated planar algebra
is isomorphic to $P^{\langle g_i : i\in I \rangle, \omega}$.
\end{corollary}
\noindent{\bf Proof:} The proof follows from
\cite{J0} and Theorem \ref{mainiso}.\qed

\begin{remark}
Note that the isomorphism $\phi$ in the proof of Theorem \ref{mainiso}
uses the $3$-cocycle $\omega$ only in the step involving the
conditional expectation onto the commutant of $M$. In particular, the
filtered $*$-algebra structure does not involve $\omega$.
\end{remark}

Analyzing the filtered $*$-algebra structure of our planar algebra,
one can easily find that the principal graph $\Gamma$ of $N \subset M$
is a Cayley-like graph. More precisely, if
$G_n = \{alt(\ul{i}):\ul{i} \in I^n\}$
for $n \geq 1$, and $G_0 =\{e\}$, then $V_n(\Gamma) = G_n \setminus
G_{n-2}$ denotes the set of vertices of $\Gamma$ which are at a
distance $n$ from the distinguished vertex for $n \geq 1$, and $V_0(\Gamma) =
\{e\}$. The number of edges between $g \in V_n(\Gamma)$ and $h \in
V_{n+1} (\Gamma)$ is $\sum_{i \in I} \delta_{g,h g_i}$ (resp. $\sum_{i
  \in I} \delta_{g g_i,h}$) if $n$ is odd (resp. even). Note that
this is well-known.

The most elegant feature of the
planar algebra $P^{\langle g_i : i\in I \rangle, \omega}$ is that
the distinguished basis forms the `loop-basis' of the filtered
$*$-algebra arising from paths on the principal graph. Note that
the $3$-cocycle $\omega$ does not enter in the definition of
the actions of multiplication, inclusion and unit tangles (defined in Section
\ref{plnalg}) or in the $*$-operation. Of course, we found the
abstract planar algebra by first computing the action of tangles
on the relative commutants. We then deduced from it an abstract
prescription of the planar algebra associated to a $G$-kernel
and a $3$-cocycle, which is the one presented in Section 3.

The path algebra associated to the principal graph can always be
used to obtain
a description of the filtered $*$-algebra structure of a subfactor
planar algebra (see for instance \cite{JS}). The extra information
encoded in the planar algebra which the principal graph
cannot provide, is the action of the left conditional expectation
tangle (or equivalently, the rotation tangle or the left-inclusion
tangle with an even number of strings). The main issue in this
paper was the choice of the unitaries $v(\ul{i})$ satisfying the
conclusion of Lemma \ref{unichoice} in such a way that the basis
elements of $N^\prime \cap M_n$  correspond to loops on the principal
graph with the correspondence extending to a filtered $*$-algebra
isomorphism. Note that this choice of $v(\ul{i})$ is unique up to
a scalar in $\mathbb T$. It is a delicate choice in the sense that
another choice would very likely make the $3$-cocycle $\omega$ appear
in the description of the filtered $*$-algebra structure, whereas
$\omega$ does not feature in the path algebra on the principal
graph $\Gamma$.

We will prove next a converse of Theorem \ref{mainiso}.
We will refer to the abstract planar algebra defined in Section \ref{plnalg} as
{\it diagonal planar algebra}.
\begin{theorem}\label{conv}
Any finite index extremal subfactor $N \subset M$ whose standard invariant is
given by a diagonal planar algebra is a diagonal subfactor. Moreover, if
the associated group and its generators of the diagonal planar algebra 
is given by $G$ and
$\{g_i:i \in I\}$ repectively, then for every $i_0 \in I$, there exists
$\alpha_i \in \Aut N$, $i \in I$, such that:

(i) $(N \subset M) \cong (N \hookrightarrow M_I \otimes N)$ where
$N \hookrightarrow M_I \otimes N$ is the diagonal subfactor with respect to
the automorphisms $\alpha_i$ for $i \in I$.

(ii) There exists a group isomorphism
$\psi : \langle \alpha_i : i \in I \rangle_{\Out N} \rightarrow \langle
g^{-1}_{i_0}g_i : i \in I \rangle \leq G$ sending $\alpha_i$ to
$g^{-1}_{i_0}g_i$ for all $i \in I$.
\end{theorem}
\noindent{\bf Proof:} Let $P$ be the planar algebra associated to
$N \subset M$, $P^\Delta$ be a diagonal planar algebra associated to $G$ and
$\{g_i:i \in I\}$, and $\phi : P^\Delta \rightarrow P$ a $*$-planar algebra
isomorphism.\\

\noindent{\it Setting up matrix units:}\\
For all $i, j \in I$ such that $g_i = g_j$, set
$q^i_j := \phi ((i,j)) \in P_1 = N^\prime \cap M$. Note that
$\underset{i \in I}{\sum} q^i_i = 1$. So, to create other off-diagonal matrix
units $q^i_j$, we partition $I = \overset{m}{\underset{n=0}{\coprod}} I_n$
such
that:

(i) $i_0 \in I_0$,

(ii) $g_i = g_j \Leftrightarrow i, j \in I_n$ for some $0 \leq n \leq m$.\\
For each $n \in \{ 1, 2, \cdots, m\}$, choose $i_n \in I^n$ and partial
isometry $q^{i_0}_{i_n} \in M$ such that
$q^{i_0}_{i_n} \left( q^{i_0}_{i_n} \right)^* = q^{i_0}_{i_0}$
and $\left(q^{i_0}_{i_n}\right)^* q^{i_0}_{i_n} = q^{i_n}_{i_n}$. (Note that
$q^{i_0}_{i_0} = \phi ((i_0, i_0))$ and $q^{i_n}_{i_n} = \phi ((i_n, i_n))$
have the same trace $\abs{I}^{-1}$). Extend $q$ by defining
$q : I \times I \rightarrow M$ by
\[I \times I \supset I^s \times I^t \ni (i,j) \overset{q}{\mapsto} q^i_j := q^
i_{
i_s} \left( q^{i_0}_{i_s} \right)^* q^{i_0}_{i_t} q^{i_t}_j = \phi((i,i_s))
\left( q^{i_0}_{i_s} \right)^* q^{i_0}_{i_t} \phi((i_t,j)) \in M\]
It is completely routine to check (using properties of partial isometry and
action of multiplication tangle in $P^\Delta$) that (i) $q$ is well-defined,
(ii) $q^i_j q^k_l = \delta_{j,k} q^i_l$, and (iii)
$\left( q^i_j \right)^* = q^j_i$.\\

\noindent{\it Finding automorphisms:}\\
Using extremality of $N \subset M$, for each $i \in I$, we get
\begin{align*}
[q^i_i M q^i_i : N q^i_i]=1 & \Rightarrow N q^i_i = q^i_i M q^i_i\\
& \Rightarrow N q^i_j \subset q^i_i M q^j_j \subset q^i_j N \text{ and } N q^
i_j \supset q^i_i M q^j_j \supset q^i_j N\\
& \Rightarrow N q^i_j = q^i_i M q^j_j = q^i_j N \Rightarrow M \cong \underset{
i,j \in I}{\bigoplus} q^i_j N
\end{align*}
where we use $[q^i_i,N]=0$ in the second implication.

For each $i \in I$, define $\alpha_i : N \rightarrow N$ by
$q^{i_0}_i x = \alpha_i (x) q^{i_0}_i$ for all $x \in N$. Since $N$ is a
factor and $[q^j_j,N]=0$ for $j \in I$, $\alpha_i$ is well-defined
and injective; surjectivity follows from
$N q^i_j = q^i_j N$. Linearity and homomorphism property of $\alpha_i$ follow
immediately, 
and we also have the identity
$x q^{i_0}_i = q^{i_0}_i \alpha^{-1}_i (x)$. To show
$q^i_j x = \left( \alpha^{-1}_i \alpha_j \right) (x) q^i_j$, note that
\begin{align*}
& q^i_{i_0} x = q^i_{i_0} x q^{i_0}_i q^i_{i_0} = q^i_i \alpha^{-1}_i (x)
q^i_{i_0} = \alpha^{-1}_i (x) q^i_{i_0}\\
\Rightarrow \; &  q^i_j x = q^i_{i_0} \alpha_j (x) q^{i_0}_j = \left(
\alpha^{-1}_i\alpha_j \right) (x) q^i_j, \text{ for all } x \in N,
\text{ for all } i, j \in I. 
\end{align*}
Using this relation, it is easy to show that $\alpha_i$ is $*$-preserving. Two
other easy consequences are $\alpha_i = \alpha_j$ if and only if $g_i = g_j$
and $\alpha_{i_0} = id_N$.\\

\noindent{\it Structure of diagonal subfactor:}\\
Define:

(i) $\kappa : \tilde{M} := M_I \otimes N \rightarrow M$ by
$\kappa(E_{i,j} \otimes x) = q^i_j \alpha^{-1}_j (x) = \alpha^{-1}_i (x)
q^i_j$,

(ii) $\lambda: M \rightarrow \tilde{M}$ by $\lambda(x) =
\underset{i,j \in I}{\sum}
E_{i,j} \otimes \lambda_{i,j} (x)$, \\
where $\lambda_{i,j} : M \rightarrow N$ is the map given by the relation
$q^i_i x q^j_j = q^i_j \alpha^{-1}_j (\lambda_{i,j}(x))$ for all $x \in M$.

Clearly, $\kappa \circ \lambda = id_M$,
$\lambda \circ \kappa = id_{\tilde{M}}$ and $\kappa$ is a $*$-isomorphism. Set
$\tilde{N} := \lambda(N) = \left\{ \left. \underset{i \in I}{\sum} E_{i,i}
\otimes
\alpha_i (x) \right| x \in N \right\} \subset \tilde{M}$.\\

This proves that $N \subset M$ is a diagonal subfactor as claimed. 
The rest of the proof pertains to \ref{conv} (ii).

\noindent{\it Matrix units for the tower of basic construction:}\\
Let $M_n$ denotes the II$_1$ factor obtained from $N \subset M$ by  
iterating the basic construction $n$ times. 
We will first define $q^{\ul{i}}_{\ul{j}} \in M_{n-1}$ for
$\ul{i}, \ul{j} \in I^n$ and $n \geq 2$ satisfying:

(i) $\left( q^{\ul{i}}_{\ul{j}} \right)^* = q^{\ul{j}}_{\ul{i}}$ for all
$\ul{i},\ul{j} \in I^n$,

(ii) $q^{\ul{i}}_{\ul{j}} \; x \; = \; alt_{\alpha} (\ul{i}, \tilde{\ul{j}})(x)
\; q^{\ul{i}}_{\ul{j}}$ for all $\ul{i},\ul{j} \in I^n$ and $x \in N$,

(iii) $q^{\ul{i}}_{\ul{i}} = \phi ((\ul{i}, \tilde{\ul{i}}))$ for all
$\ul{i} \in I^n$,

(iv)
$q^{\ul{i}}_{\ul{j}} \; q^{\ul{k}}_{\ul{l}} \; = \; \delta_{\ul{j},\ul{k}}
\; q^{\ul{i}}_{\ul{l}}$ for all $\ul{i},\ul{j}, \ul{k}, \ul{l} \in I^n$,

(v) $q^{(\ul{s},s)}_{(\ul{t},t)} \; q^{\ul{u}}_{\ul{v}} \; = \;
\delta_{\ul{t},\ul{u}} \; q^{(\ul{s},s)}_{(\ul{v},t)}$ for all
$\ul{s},\ul{t}, \ul{u}, \ul{v} \in I^{n-1}$ and $s,t \in I$\\
by induction on $n$ where $alt_{\alpha}$ is defined by
$alt_{\alpha} (i_1, \cdots, i_m) = \alpha^{-1}_{i_1} \alpha_{i_2}
\alpha^{-1}_{i_3} \cdots \alpha^{(-1)^m}_{i_m}$.

Suppose we have defined such $q^{\ul{i}}_{\ul{j}} \;$'s
for all
$\ul{i}, \ul{j} \in I^m$ and $m \leq n$. Now, for $i,j \in I$ and
$\ul{i}, \ul{j} \in I^n$,
set $q^{(\ul{i},i)}_{(\ul{j},j)} := q^{\ul{i}}_{(\ul{s},i)} \; E_n \;
q^{(\ul{s},j)}_{\ul{j}} \in M_n$ for some $\ul{s} \in I^{n-1}$, where 
$E_n = \abs{I} e_n$ is the element in $M_n$ given by the $n$-th Jones
projection tangle. To show that
the definition of  $q^{(\ul{i},i)}_{(\ul{j},j)}$ is independent of $\ul{s} \in
I^{n-1}$, observe that
\[
q^{(\ul{i},i)}_{(\ul{j},j)} = q^{\ul{i}}_{(\ul{s},i)} \;
\phi(\ul{s}, \tilde{\ul{s}}) \;
E_n \; q^{(\ul{s},j)}_{\ul{j}} =  q^{\ul{i}}_{(\ul{s},i)} \;
q^{\ul{s}}_{\ul{t}}
\; E_n \; q^{\ul{t}}_{\ul{s}} \; q^{(\ul{s},j)}_{\ul{j}} =
q^{\ul{i}}_{(\ul{t},i)} \; E_n \; q^{(\ul{t},j)}_{\ul{j}}
\]
for all $\ul{t} \in I^{n-1}$. Properties (i) and (v) hold trivially. For (ii),
note that
\[
q^{(\ul{i},i)}_{(\ul{j},j)} \; x = alt_{\alpha} (\ul{i}, i, \tilde{\ul{s}},
\ul{s},j, \tilde{\ul{j}})(x)
\; q^{(\ul{i},i)}_{(\ul{j},j)} = alt_{\alpha} ((\ul{i}, i),
\widetilde{(\ul{j},j)})(x) \; q^{(\ul{i},i)}_{(\ul{j},j)} \text{ for all } x
\in N.
\]

Next, we prove property (iv). For $i \in I$ and $m \geq 1$, let $\eta_m (i)$
denote the element
\psfrag{eta}{$\eta_m (i)$ =}
\psfrag{2}{$2$}
\psfrag{m-1}{$m-1$}
\psfrag{m}{$m$}
\psfrag{}{$$}
\psfrag{qii}{$q^i_i$}
\psfrag{cond1}{if $m$ is odd}
\psfrag{cond2}{if $m$ is even}
\begin{center}
\epsfig{file=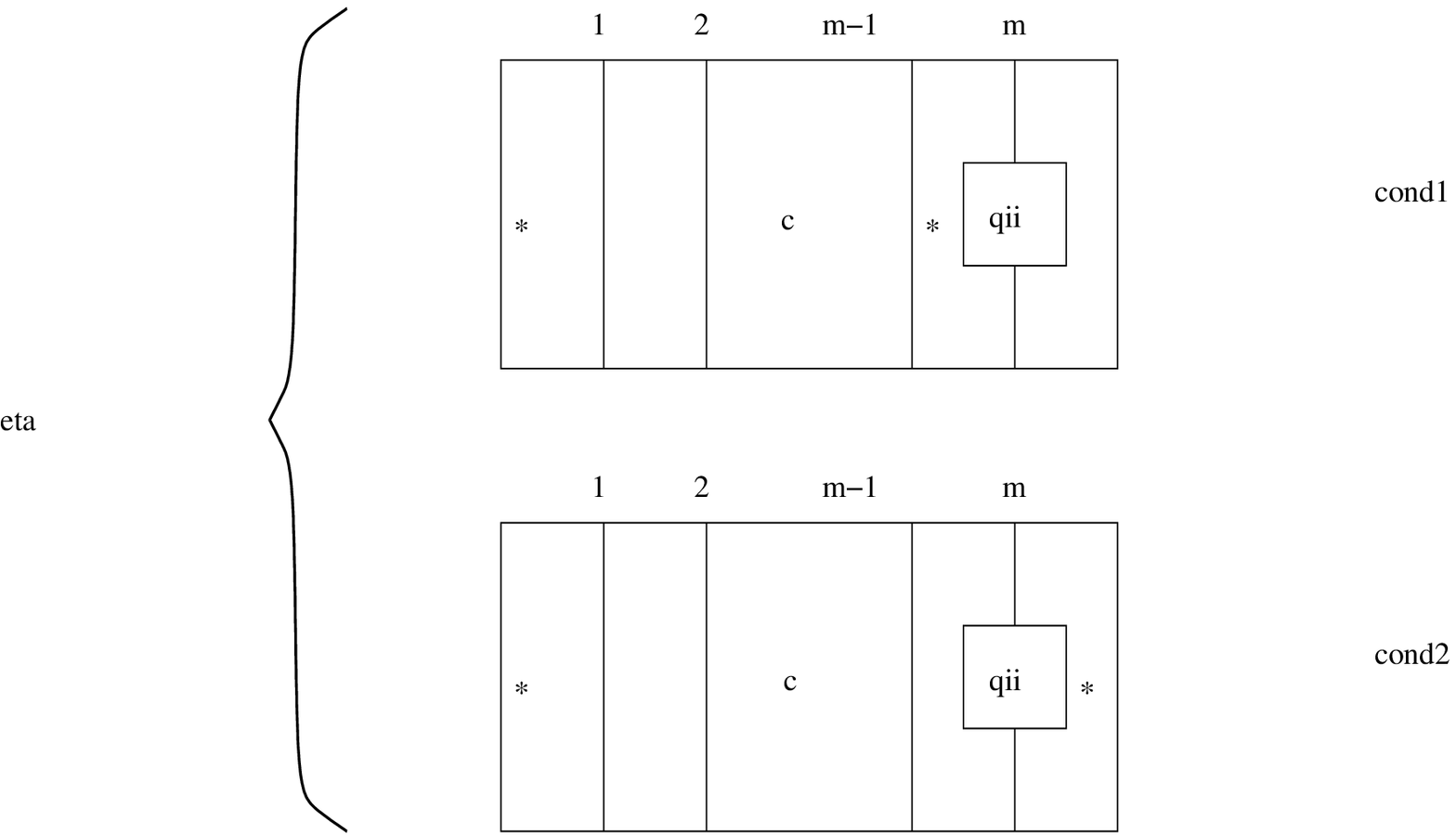, scale=0.5}
\end{center}
in $M^{\prime}_{m-2} \cap M_{m-1}$ 
~(where $M_{-1} := N$, $M_0 := M$).
Two important relations which we will
often use, are $\eta_m (i) E_m = \eta_{m+1} (i) E_m$ and $E_m \eta_m (i) =
E_m \eta_{m+1} (i)$. Getting back to property (iv), we have
\begin{align*}
& q^{(\ul{i},i)}_{(\ul{j},j)} \; q^{(\ul{k},k)}_{(\ul{l},l)} =
q^{\ul{i}}_{(\ul{s},i)} \; E_n
\; q^{(\ul{s},j)}_{\ul{j}} \; q^{\ul{k}}_{(\ul{s},k)} \; E_n \;
q^{(\ul{s},l)}_{\ul{l}}
= \delta_{\ul{j},\ul{k}} \;
q^{\ul{i}}_{(\ul{s},i)} \; E_n
\; q^{(\ul{s},j)}_{(\ul{s},k)} \; E_n \;
q^{(\ul{s},l)}_{\ul{l}}\\
= \; & \delta_{\ul{j},\ul{k}} \;
q^{\ul{i}}_{(\ul{s},i)} \; E_n \; \eta_n (j)
\; q^{(\ul{s},j)}_{(\ul{s},k)} \; E_n
\; q^{(\ul{s},l)}_{\ul{l}}
=  \delta_{\ul{j},\ul{k}} \;
q^{\ul{i}}_{(\ul{s},i)} \; E_n \; \eta_{n+1} (j)
\; q^{(\ul{s},j)}_{(\ul{s},k)} \; E_n
\; q^{(\ul{s},l)}_{\ul{l}}\\
= \; & \delta_{\ul{j},\ul{k}} \;
q^{\ul{i}}_{(\ul{s},i)} \; E_n
\; q^{(\ul{s},j)}_{(\ul{s},k)} \; \eta_{n+1} (j) \; E_n
\; q^{(\ul{s},l)}_{\ul{l}}
= \delta_{\ul{j},\ul{k}} \;
q^{\ul{i}}_{(\ul{s},i)} \; E_n
\; q^{(\ul{s},j)}_{(\ul{s},k)} \; \eta_{n} (j) \; E_n
\; q^{(\ul{s},l)}_{\ul{l}}\\
= \; & \delta_{\ul{j},\ul{k}} \; \delta_{j,k} \;
q^{\ul{i}}_{(\ul{s},i)} \; E_n
\; q^{(\ul{s},j)}_{(\ul{s},j)} \; E_n
\; q^{(\ul{s},l)}_{\ul{l}}
= \delta_{\ul{j},\ul{k}} \; \delta_{j,k} \;
q^{\ul{i}}_{(\ul{s},i)} \; E_n \; \phi((\ul{s}, j, j , \tilde{\ul{s}})) \; E_n
\; q^{(\ul{s},l)}_{\ul{l}}\\
= \; & \delta_{\ul{j},\ul{k}} \; \delta_{j,k} \;
q^{\ul{i}}_{(\ul{s},i)} \; \phi((\ul{s},\tilde{\ul{s}})) \; E_n
\; q^{(\ul{s},l)}_{\ul{l}}
= \delta_{\ul{j},\ul{k}} \; \delta_{j,k} \;
q^{\ul{i}}_{(\ul{s},i)} \; E_n
\; q^{(\ul{s},l)}_{\ul{l}}
= \delta_{\ul{j},\ul{k}} \; \delta_{j,k} \;
q^{(\ul{i},i)}_{(\ul{l},l)}.
\end{align*}

It remains to establish property (iii). Now, for $\ul{i} \in I^{n}$, $i \in I$
and $\ul{s} \in I^{n-1}$,
\begin{align*}
& \phi ((\ul{i}, i, i, \tilde{\ul{i}})) \; q^{(\ul{i},i)}_{(\ul{i},i)} \;
\phi ((\ul{i}, i, i, \tilde{\ul{i}}))
= \phi ((\ul{i}, \tilde{\ul{i}})) \; \eta_{n+1} (i) \; q^{\ul{i}}_{\ul{s},i}
\; E_n \; q^{\ul{s},i}_{\ul{i}} \; \eta_{n+1} (i) \;
\phi ((\ul{i}, \tilde{\ul{i}}))\\
= \; & \phi ((\ul{i}, \tilde{\ul{i}})) \;  q^{\ul{i}}_{\ul{s},i} \;
\eta_{n+1} (i) \; E_n \; \eta_{n+1} (i) \; q^{\ul{s},i}_{\ul{i}} \;
\phi ((\ul{i}, \tilde{\ul{i}}))
= q^{\ul{i}}_{\ul{s},i} \;
\eta_n (i) \; E_n \; \eta_n (i) \; q^{\ul{s},i}_{\ul{i}}\\
= \; & q^{\ul{i}}_{\ul{s},i} \; E_n \; q^{\ul{s},i}_{\ul{i}}
= q^{(\ul{i},i)}_{(\ul{i},i)}.
\end{align*}
Since $q^{(\ul{i},i)}_{(\ul{i},i)} \in \mathcal{P} (N^\prime \cap M_n)$ (by
property (ii)) and
$\phi ((\ul{i}, i, i, \tilde{\ul{i}}))$ is a minimal projection of
$P_{n+1} = N^\prime \cap M_n$, therefore
$q^{(\ul{i},i)}_{(\ul{i},i)} = \phi ((\ul{i}, i, i, \tilde{\ul{i}}))$.

The proof for the initial case of $n = 2$ is similar and is left to
reader.\\

\noindent{\it Correspondence between relations satisfied by $\alpha_i$'s and
$g_i$'s:}\\
In this part, we will prove that for $\ul{i},\ul{j} \in I^n$,
$alt (\ul{i},\tilde{\ul{j}}) = e$ if and only if $alt_\alpha
(\ul{i},\tilde{\ul{j}}) \in \Int N$.

Following the construction of the isomorphism $\lambda$ between $M$ and
$\tilde{M}$, one can define an isomorphism $\lambda^{(n)}$ as follows:
\[
M_{n-1} \ni x \overset{\lambda^{(n)}}{\longmapsto} \lambda^{(n)} (x) =
\underset{\ul{k},\ul{l} \in I^n}{\sum} E_{\ul{k},\ul{l}} \otimes
\lambda^{(n)}_{\ul{k},\ul{l}} (x) \in \tilde{M}_{n-1} := M_{I^n} \otimes N
\]
where $\lambda^{(n)}_{\ul{k},\ul{l}} : M_{n-1} \rightarrow N$ is the map given
by the relation
\[
q^{\ul{k}}_{\ul{k}} \; x \; q^{\ul{l}}_{\ul{l}} = q^{\ul{k}}_{\ul{l}} \;
alt_\alpha (\ul{l}) \left( \lambda^{(n)}_{\ul{k},\ul{l}} (x) \right).
\]
Thus, $\lambda^{(n)} (N) = \left\{ \underset{\ul{k} \in I^n}{\sum}
E_{\ul{k},\ul{k}} \otimes alt_\alpha^{-1} (\ul{k}) (x)  \; : \; x \in
N \right\}$, for $n \geq 1$. Note that $\lambda^{(1)} = \lambda$.

Let $alt (\ul{i},\tilde{\ul{j}}) = e$ for $\ul{i},\ul{j} \in I^n$. Note that
$\lambda^{(n)} (q^{\ul{i}}_{\ul{j}}) = E_{\ul{i}, \ul{j}} \otimes 1$ and
$\lambda^{(n)} \left(\phi((\ul{i}, \tilde{\ul{j}})) \right)$ are partial
isometries between
$\lambda^{(n)} (q^{\ul{i}}_{\ul{i}}) = \lambda^{(n)} \left(\phi((\ul{i},
\tilde{\ul{i}})) \right) =  E_{\ul{i}, \ul{i}} \otimes 1$ and
$\lambda^{(n)} (q^{\ul{j}}_{\ul{j}}) = \lambda^{(n)} \left(\phi((\ul{j},
\tilde{\ul{j}})) \right) = E_{\ul{j}, \ul{j}} \otimes 1$. So, there exists
$u \in \mathcal{U} (N)$ such that
$\lambda^{(n)} \left(\phi((\ul{i}, \tilde{\ul{j}})) \right) = E_{\ul{i},
\ul{j}} \otimes u =
\lambda^{(n)} \left( q^{\ul{i}}_{\ul{j}} \; v \right) $ where
$v = alt_\alpha (\ul{j}) (u) \in \mathcal{U} (N)$. Hence,
\begin{align*}
& q^{\ul{i}}_{\ul{j}} \; v = \phi((\ul{i}, \tilde{\ul{j}})) \in N^\prime \cap
M_{n-1}\\
\Rightarrow \; & y \; q^{\ul{i}}_{\ul{j}} \; v = q^{\ul{i}}_{\ul{j}} \;
v \; y = \left( alt_\alpha ((\ul{i}, \tilde{\ul{j}})) \circ \Ad v \right) (y)
\; q^{\ul{i}}_{\ul{j}} \; v \; \text{ for all } y \in N\\
\Rightarrow \; & \left( alt_\alpha ((\ul{i}, \tilde{\ul{j}})) \circ
\Ad v \right) (y)
= y \; \text{ for all } y \in N\\
\Rightarrow \; & alt_\alpha ((\ul{i}, \tilde{\ul{j}})) \in \Int N.
\end{align*}

Conversely, if $alt_\alpha ((\ul{i}, \tilde{\ul{j}})) \in \Int N$, that is,
$alt_\alpha ((\ul{i}, \tilde{\ul{j}})) \circ \Ad v = id_N$ for
some $v \in \mathcal{U} (N)$, then $\left(  E_{\ul{i},\ul{j}} \otimes
alt^{-1}_\alpha (\ul{j}) (v) \right) \in \left( \left( \lambda^{(n)} (N)
\right)^\prime \cap \tilde{M}_{n-1} \right)$.
Now, $alt ((\ul{i}, \tilde{\ul{j}})) \neq e$ implies $dim
\left( \left( \lambda^{(n)} (N) \right)^\prime \cap
\tilde{M}_{n-1} \right)
> dim \left( N^\prime \cap M_{n-1} \right) $ which is a contradiction. Hence,
$alt ((\ul{i}, \tilde{\ul{j}})) = e$.\\

\noindent{\it The group generated by $\alpha_i$'s:}\\
Let $H := \langle \theta_i = [\alpha_i]_{\Out N} : i \in I \rangle \leq
{\Out N}$,
$\tilde{H} := \langle g^{-1}_{i_0} g_i : i \in I \rangle \leq {G}$, 
$J := I \times \{1,-1\}$.  
Define maps
$w_H : \underset{n \geq 1}{\coprod} J^n \rightarrow H$ and 
$w_{\tilde{H}} :\underset{n \geq 1}{\coprod} J^n \rightarrow \tilde{H}$  
by 
\begin{align*}
J^n \ni ((i_1, \epsilon_1), (i_2, \epsilon_2), \cdots (i_n,
\epsilon_n)) & \overset{w_H}{\longmapsto} \theta_{i_1}^{\epsilon_1}
\theta_{i_2}^{\epsilon_2} \cdots \theta_{i_n}^{\epsilon_n} \in H\\
J^n \ni ((i_1, \epsilon_1), (i_2, \epsilon_2), \cdots (i_n,
\epsilon_n)) & \overset{w_{\tilde{H}}}{\longmapsto} 
({g_{i_0}^{-1}g_{i_1}})^{\epsilon_1}
({g_{i_0}^{-1}g_{i_2}})^{\epsilon_2}
\cdots ({g_{i_0}^{-1}g_{i_n}})^{\epsilon_n} \in \tilde{H}
\end{align*}
where $\epsilon_i \in \{1,-1\}$ for $1 \leq i \leq n$.
Define $\gamma: H \rightarrow \tilde{H}$ by
$\gamma \left( w_H (\ul{j}) \right) = w_{\tilde{H}} (\ul{j})$. 
For $\gamma$ to be an isomorphism, it is enough to show $\gamma$ is well-defined
and injective. Suppose the map $\rho : J \rightarrow I^2$ sends $(i,1)$ (resp.
$(i,-1)$) to $(i_0,i)$ (resp. $(i, i_0)$). Extend $\rho$ to
$\rho : J^n \rightarrow I^{2n}$ entrywise. Note that
$w_H (\ul{j}) = alt_H \left( \rho(\ul{j}) \right)$ and
$w_{\tilde{H}} (\ul{j}) = alt \left( \rho(\ul{j}) \right)$. This implies
\[
w_H (\ul{j}) = 1_H \Leftrightarrow alt_\alpha \left( \rho(\ul{j}) \right)
\in \Int N \Leftrightarrow alt \left( \rho(\ul{j}) \right) = e
\Leftrightarrow w_{\tilde{H}} (\ul{j}) = e.
\]
Hence, $H$ and $\tilde{H}$ are isomorphic.
\qed
\bibliographystyle{amsplain}

\end{document}